\documentclass[12pt,twoside,a4paper]{article}
\usepackage{bm}
\usepackage{braket}
\usepackage{graphicx}
\usepackage{graphics}
\usepackage{theorem}
\usepackage{amsmath}
\usepackage{amssymb,mathrsfs}
\usepackage{latexsym}
\usepackage{cases}
\usepackage[dvipdfmx,bookmarksnumbered,colorlinks,linkcolor=blue,urlcolor=blue,citecolor=blue]{hyperref}
\usepackage{color}	
\usepackage{algorithm}
\usepackage{algpseudocode}
\usepackage{caption}
\newfloat{procedure}{H}{loa}
\floatname{procedure}{Procedure} 
\newtheorem{THM}{Theorem}[section]
\newtheorem{LEM}[THM]{Lemma}
\newtheorem{PROP}[THM]{Proposition}
\newtheorem{COR}[THM]{Corollary}
\newtheorem{CON}{Condition}
\newtheorem{ASSU}[THM]{Assumption}

\newtheorem{NOTI}[THM]{Notion}

\theorembodyfont{\rmfamily}
\newtheorem{REM}[THM]{Remark}
\newtheorem{DEF}[THM]{Definition}

\newtheorem{PROB}[THM]{Problem}
\newenvironment{proof}{\noindent{\it Proof.~~}}{\qed}
\newcommand{\qed}{\hspace*{\fill}$\Box$}

\numberwithin{equation}{section}  
\newcommand{\vc}{\bm}
\makeatletter
\DeclareRobustCommand\widecheck[1]{{\mathpalette\@widecheck{#1}}}
\def\@widecheck#1#2{%
    \setbox\z@\hbox{\m@th$#1#2$}%
    \setbox\tw@\hbox{\m@th$#1%
       \widehat{%
          \vrule\@width\z@\@height\ht\z@
          \vrule\@height\z@\@width\wd\z@}$}%
    \dp\tw@-\ht\z@
    \@tempdima\ht\z@ \advance\@tempdima2\ht\tw@ \divide\@tempdima\thr@@
    \setbox\tw@\hbox{%
       \raise\@tempdima\hbox{\scalebox{1}[-1]{\lower\@tempdima\box
\tw@}}}%
    {\ooalign{\box\tw@ \cr \box\z@}}}
\makeatother
\makeatletter
\newif\if@borderstar
\def\bordermatrix{\@ifnextchar*{%
 \@borderstartrue\@bordermatrix@i}{\@borderstarfalse\@bordermatrix@i*}%
}
\def\@bordermatrix@i*{\@ifnextchar[{\@bordermatrix@ii}{\@bordermatrix@ii[()]}}
\def\@bordermatrix@ii[#1]#2{%
\begingroup
 \m@th\@tempdima8.75\p@\setbox\z@\vbox{%
 \def\cr{\crcr\noalign{\kern 2\p@\global\let\cr\endline }}%
 \ialign {$##$\hfil\kern 2\p@\kern\@tempdima & \thinspace %
  \hfil $##$\hfil && \quad\hfil $##$\hfil\crcr\omit\strut %
  \hfil\crcr\noalign{\kern -\baselineskip}#2\crcr\omit %
  \strut\cr}}%
 \setbox\tw@\vbox{\unvcopy\z@\global\setbox\@ne\lastbox}%
 \setbox\tw@\hbox{\unhbox\@ne\unskip\global\setbox\@ne\lastbox}%
 \setbox\tw@\hbox{%
  $\kern\wd\@ne\kern -\@tempdima\left\@firstoftwo#1%
  \if@borderstar\kern 2pt\else\kern -\wd\@ne\fi%
 \global\setbox\@ne\vbox{\box\@ne\if@borderstar\else\kern 2\p@\fi}%
 \vcenter{\if@borderstar\else\kern -\ht\@ne\fi%
  \unvbox\z@\kern -\if@borderstar2\fi\baselineskip}%
 \if@borderstar\kern-2\@tempdima\kern2\p@\else\,\fi\right\@secondoftwo#1 $%
 }\null \;\vbox{\kern\ht\@ne\box\tw@}%
\endgroup
}
\makeatother

\newcommand{\ol}{\overline}

\newcommand{\wt}{\widetilde}
\newcommand{\wh}{\widehat}
\newcommand{\bv}{\breve}

\newcommand{\down}[2]{\smash{\lower#1\hbox{#2}}}
\newcommand{\up}[2]{\smash{\lower-#1\hbox{#2}}}

\newcommand{\dm}{\displaystyle}

\def\presub#1{\hspace{0.05em}{}_{#1}\hspace{-0.05em}}

\newcommand{\EE}{\mathbb{E}}
\newcommand{\PP}{\mathbb{P}}

\newcommand{\bbA}{\mathbb{A}}

\newcommand{\bbC}{\mathbb{C}}

\newcommand{\bbG}{\mathbb{G}}
\newcommand{\bbH}{\mathbb{H}}
\newcommand{\bbI}{\mathbb{I}}
\newcommand{\bbJ}{\mathbb{J}}
\newcommand{\bbK}{\mathbb{K}}
\newcommand{\bbO}{\mathbb{O}}
\newcommand{\bbL}{\mathbb{L}}
\newcommand{\bbM}{\mathbb{M}}
\newcommand{\bbN}{\mathbb{N}}

\newcommand{\bbS}{\mathbb{S}}

\newcommand{\bbX}{\mathbb{X}}

\newcommand{\bbZ}{\mathbb{Z}}
\DeclareMathOperator*{\argmin}{arg\,min}
\DeclareMathOperator*{\argmax}{arg\,max}

\newcommand{\diag}{\mathrm{diag}}

\newcommand{\row}{\mathrm{row}}

\newcommand{\varep}{\varepsilon}

\renewcommand{\labelenumi}{(\roman{enumi})}
\newcommand{\dd}[1]{\if#11 1\!\!1 
\else {\if#1C I\!\!\!C
\else {\if#1G I\!\!\!G 
\else {\if#1J J\!\!\!J 
\else {\if#1S S\!\!\!S
\else {\if#1Z Z\!\!\!Z
\else {\if#1Q O\!\!\!\!Q
\else I\!\!#1
\fi}
\fi}
\fi}
\fi} 
\fi} 
\fi} 
\fi} 
\makeatletter
\def\eqnarray{\stepcounter{equation}\let\@currentlabel=\theequation
\global\@eqnswtrue
\global\@eqcnt\z@\tabskip\@centering\let\\=\@eqncr
$$\halign to \displaywidth\bgroup\@eqnsel\hskip\@centering
  $\displaystyle\tabskip\z@{##}$&\global\@eqcnt\@ne 
  \hfil$\;{##}\;$\hfil
  &\global\@eqcnt\tw@ $\displaystyle\tabskip\z@{##}$\hfil 
   \tabskip\@centering&\llap{##}\tabskip\z@\cr}
\makeatother

\makeatother

\begin{document}\thispagestyle{empty} 

\hfill

\vspace{-10mm}

{\large{\bf
\begin{center}
A new matrix-infinite-product-form solution for upper block-Hessenberg Markov chains and its quasi-algorithmic constructibility%
\footnote[1]{%
To appear in Advances in Applied Probability, vol. 55, no. 1, March 2023
}
%
%
\end{center}
}
}

\begin{center}
{
Hiroyuki Masuyama%
\footnote[2]{E-mail: masuyama@sys.i.kyoto-u.ac.jp}
}

\medskip

{\small
Graduate School of Management, Tokyo Metropolitan University, Tokyo 192--0364, Japan.
}

\bigskip
\medskip

{\small
\textbf{Abstract}

\medskip

\begin{tabular}{p{0.85\textwidth}}
This paper presents a new matrix-infinite-product-form (MIP-form) solution for the stationary distribution in upper block-Hessenberg Markov chains (UBH-MCs). The existing MIP-form solution (Masuyama, Queueing Syst., Vol.~92, 2019, pp.~173--200) requires a certain parameter set that satisfies both a Foster-Lyapunov drift condition and a convergence condition. In contrast, the new MIP-form solution requires no such parameter sets and no other conditions. The new MIP-form solution also has ``quasi-algorithmic constructibility", which is a newly introduced feature of being constructed by iterating infinitely many times a recursive procedure of finite complexity per iteration. This feature is not found in the other solutions for the stationary distribution in UBH-MCs.
\end{tabular}
}
\end{center}

\begin{center}
\begin{tabular}{p{0.90\textwidth}}
{\small

\medskip

}
\end{tabular}

\end{center}

\section{Introduction}\label{sec-intro}

This paper studies the stationary distribution in continuous-time upper block-Hessenberg Markov chains (UBH-MCs). This study focuses on establishing a {\it theoretical} procedure for constructing the {\it exact} stationary distribution in an algorithmic way, rather than establishing a {\it practical} procedure for computing an {\it approximate} stationary distribution (of course, such practicality is a significant issue, too). Additionally, this study would be a {\it partial} answer to the question: ``What is a class of Markov chains whose {\it exact} stationary distributions are theoretically constructible by an algorithmic procedure?"

The class of UBH-MCs is characterized by the upper block-Hessenberg (UBH) structure of the infinitesimal generator matrix (called {\it generator} for short), and this class includes M/G/1-type Markov chains and level-dependent quasi-birth-and-death processes (LD-QBDs). Remarkably, even a multi-dimensional random walk on the nonnegative lattice is expressed as a UBH-MC such that its level variable is the sum of coordinate components and its phase variable is the set of all the coordinate components but one. Thus, the class of UBH-MCs is a powerful tool for modeling information and communication systems, inventory systems, transportation systems, etc. Moreover, the time-average performance of these systems is evaluated through the stationary distribution. That is why this study focuses on the stationary distribution in UBH-MCs.

We now introduce the definition of UBH-MCs in order to describe the background and purpose of this study. Let $\bbS$ denote a countable set such that
\[
\bbS = \bigcup_{k \in \bbZ_+} \{k\} \times \bbM_k,
\]
where $\bbZ_+=\{0,1,2,\dots\}$ and $\bbM_k = \{1,2,\dots,M_k\} \subset \bbN:=\{1,2,3,\dots\}$. Let $\vc{Q}:=(q(k,i;\ell,j))_{(k,i;\ell,j) \in\bbS^2}$ denote
an essential $Q$-matrix (see Definition~\ref{defn-essential}), and thus
\begin{alignat*}{2}
-\infty < q(k,i;k,i) &< 0, &\qquad  (k,i) &\in \bbS,
\\ 
0 \le q(k,i;\ell,j) &< \infty, &\qquad  (k,i) &\in \bbS,~(\ell,j) \in \bbS \setminus\{(k,i)\},
\\
\sum_{(\ell,j) \in \bbS} q(k,i;\ell,j) &= 0, & (k,i) &\in \bbS.
\end{alignat*}
Assume that $\vc{Q}$ is of UBH form:
\begin{equation}
\vc{Q} = 
\bordermatrix{
        & \bbL_0 & \bbL_1    & \bbL_2    & \bbL_3   & \cdots
\cr
\bbL_0 		& 
\vc{Q}_{0,0}  	& 
\vc{Q}_{0,1} 	& 
\vc{Q}_{0,2}  	& 
\vc{Q}_{0,3} 	&  
\cdots
\cr
\bbL_1 		&
\vc{Q}_{1,0}  	& 
\vc{Q}_{1,1} 	& 
\vc{Q}_{1,2}  	& 
\vc{Q}_{1,3} 	&  
\cdots
\cr
\bbL_2 		& 
\vc{O}  		& 
\vc{Q}_{2,1} 	& 
\vc{Q}_{2,2}  	& 
\vc{Q}_{2,3} 	&  
\cdots
\cr
\bbL_3 & 
\vc{O}  		& 
\vc{O}  		& 
\vc{Q}_{3,2}  	& 
\vc{Q}_{3,3} 	&  
\cdots
\cr
~\vdots  	& 
\vdots     		& 
\vdots     		&  
\vdots    		& 
\vdots    		& 
\ddots
},
\label{defn-Q}
\end{equation}
where $\bbL_k := \{k\} \times \bbM_k$ is called {\it level} $k$ and an element $i$ of $\bbM_k$ is called {\it phase $i$ (of level $k$}). 
A Markov chain having this $Q$-matrix $\vc{Q}$ as its generator is said to be an upper block-Hessenberg Markov chain (UBH-MC). Note that a UBH-MC may be called a {\it level-dependent M/G/1-type Markov chain}. Note also that if $\vc{Q}_{k,k+m} = \vc{O}$ for all $m \ge 2$ and $k \in \bbZ_+$ then the generator is of block-tridiagonal form and thus the UBH-MC is reduced to an LD-QBD. 

We provide a basic assumption and some definitions associated with the stationary distribution vector of the UBH generator $\vc{Q}$ given in (\ref{defn-Q}). Assume that $\vc{Q}$ is ergodic (i.e., irreducible and positive recurrent) throughout the paper, unless otherwise stated. Let $\vc{\pi}:=(\pi(k,i))_{(k,i)\in\bbS}$ denote the unique and positive stationary distribution vector of the ergodic generator $\vc{Q}$ (see, e.g., \cite[Chapter 5, Theorems 4.4 and 4.5]{Ande91}). By definition,
\[
\vc{\pi}\vc{Q} = \vc{0}, \quad \vc{\pi}\vc{e}=1, \quad \vc{\pi}>\vc{0},
\]
where $\vc{e}$ denotes a column vector of ones with an appropriate dimension. For later reference, let $\vc{\pi}_k=(\pi(k,i))_{i\in\bbM_k}$ for $k \in \bbZ_+$, which leads to the level-wise partition of~$\vc{\pi}$:
\begin{eqnarray*}
\vc{\pi}
&=& 
\bordermatrix{
       	& 
\bbL_0 	& 
\bbL_1 	& 
\cdots
\cr
 	& 
\vc{\pi}_0 & 
\vc{\pi}_1 & 
\cdots
}.
\end{eqnarray*}
Furthermore, let $\vc{\pi}^{(n)}$, $n \in \bbZ_+$, denote
\begin{equation}
\vc{\pi}^{(n)} 
= {
(\vc{\pi}_0,\vc{\pi}_1,\dots,\vc{\pi}_n) \over \sum_{\ell=0}^n \vc{\pi}_{\ell}\vc{e}
},
\label{defn-pi^{(N)}}
\end{equation}
which is referred to as the {\it (finite-level) conditional stationary distribution vector}. For any sufficiently large $n$, $\vc{\pi}^{(n)}$ can be considered an approximation to $\vc{\pi}$.

The main purpose of this paper is to present a {\it quasi-algorithmically constructible} solution for the stationary distribution vector $\vc{\pi}$ of the ergodic UBH generator $\vc{Q}$ given in (\ref{defn-Q}). We here provide our notion (not strict definition) of {\it quasi-algorithmic constructions} for the stationary distribution (the stationary distribution vector) in ergodic countable-state Markov chains, which are not, of course, restricted to UBH-MCs.
\begin{NOTI}\label{defn-construction}
({\bf Quasi-algorithmic constructions for the stationary distribution in countable-state Markov chains}.) Consider an ergodic countable-state Markov chain with generator (essential $Q$-matrix) $\cal{Q}$. Let $\pi$ denote the stationary distribution of the ergodic Markov chain. Let $\mathscr{P}$ denote a procedure for sequentially generating tentative solutions $\pi(0),\pi(1),\pi(2),\dots$ for the stationary distribution $\pi$. The procedure $\mathscr{P}$ is said to be a {\it quasi-algorithmic construction} of $\pi$ and furthermore $\pi$ is said to be {\it quasi-algorithmic constructible (or have quasi-algorithmic constructibility)}, if $\mathscr{P}$ has the following features:
\begin{enumerate}
\item The procedure $\mathscr{P}$ generates each tentative solution $\pi(n)$ by using at most a finite (not necessarily bounded) number of elements of $\cal{Q}$ together with (if required) some or all the components of the previous tentative solutions $\{\pi(m);m \in \bbZ_{[0,n-1]}\}$, where $\bbZ_{[0,k]} = \{0,1,\dots,k\}$ for $k \in \bbZ_+$.
\item It takes at most {\it finite complexity} for $\mathscr{P}$ to generate each tentative solution $\pi(n)$.

\item The sequence $\{\pi(n);n\in\bbZ_+\}$ converges to $\pi$ in some mathematical sense (e.g., in the $\ell_1$-distance).
\end{enumerate}
\end{NOTI}

\begin{REM}
To the best of our knowledge, there is no universal formal definition of ``algorithm". Indeed, the term of ``algorithm" is widely used with some ambiguity in various contexts. Nevertheless, it would be generally accepted, as an informal definition, that an ``algorithm" is a finite set of operations (of finite computational complexity) to accomplish a particular task. Considering this situation, we introduce the notion of quasi-algorithmic constructions in order to distinguish our new solution presented in this paper from the existing ones for the stationary distribution in UBH-MCs.
\end{REM}

\begin{REM}\label{rem-QAC}
Any quasi-algorithmic construction can be implemented as an ``algorithm" if the construction is equipped with an appropriate stopping criterion so that it stops after finitely many iterations, and if it is possible to store the elements required to compute each tentative solution.
\end{REM} 

As far as we know, there have been no studies that achieve a quasi-algorithmic construction of the stationary distribution vector in UBH-MCs including LD-QBDs. There are, however, many related studies on the computation of the stationary distribution vector in these Markov chains.

Bright and Taylor~\cite{Brig95} presented a nice {\it matrix-product form} of 
the stationary distribution vector $\vc{\pi}=(\vc{\pi}_0,\vc{\pi}_1,\dots)$ of the LD-QBD generator:
\begin{align}
\left\{
\begin{array}{l}
\vc{\pi}_k = \vc{\pi}_0 \vc{R}_1 \vc{R}_2 \cdots  \vc{R}_k,
\qquad k \in \bbN,
\\
\vc{\pi}_0 \left(\vc{Q}_{0,0} + \vc{R}_1 \vc{Q}_{1,0}\right) = \vc{0},
\\
\vc{\pi}_0\vc{e}
+ 
\vc{\pi}_0 \dm\sum_{k=1}^{\infty}\vc{R}_1 \vc{R}_2 \cdots  \vc{R}_k\vc{e}=1.
\end{array}
\right.
\label{matrix-product-form-LD-QBD}
\end{align}
The set of matrices $\{\vc{R}_k;k\in\bbN\}$ is the minimal nonnegative solution for the system of matrix equations:
\begin{align}
\vc{Q}_{k-1,k} + \vc{R}_k \vc{Q}_{k,k} + \vc{R}_k\vc{R}_{k+1} \vc{Q}_{k+1,k}
= \vc{O},\qquad k \in \bbN.
\label{recursion-R_k}
\end{align}

The matrix-product form (\ref{matrix-product-form-LD-QBD}) appears easy to compute. Thus, based on it, Bright and Taylor~\cite{Brig95} proposed a procedure for approximately computing the stationary distribution vector of the LD-QBD generator. Besides, some researchers use the matrix-product form (\ref{matrix-product-form-LD-QBD}) as the foundation of their respective algorithms (see, e.g., \cite{Baum12-Procedia,Phun10-QTNA}).

However, the matrix-product form (\ref{matrix-product-form-LD-QBD}) does not lead to any quasi-algorithmic construction of the stationary distribution vector in LD-QBDs. Equation (\ref{recursion-R_k}) yields
\begin{align*}
\vc{R}_k 
= \vc{Q}_{k-1,k} (-\vc{Q}_{k,k} - \vc{R}_{k+1} \vc{Q}_{k+1,k})^{-1},\qquad k \in \bbN.
\end{align*}
This recursive formula shows that each $\vc{R}_k$ can be computed provided that $\vc{R}_{k+1}$ is given. Therefore, to implement the algorithms based on (\ref{matrix-product-form-LD-QBD}), we have to truncate the infinite sequence $\{\vc{R}_k; k \in \bbZ_+\}$ at a sufficiently large $K^* \in \bbZ_+$ by letting $\vc{R}_{K^*+1} = \vc{O}$. Moreover, owing to this implementation, computing {\it better} approximations (to the stationary distribution) requires setting $R$-matrices of {\it higher} levels to be zero matrices, which implies that all the components of such an approximation are computed anew from scratch.

Several researchers have studied the approximate computation of the stationary distribution in UBH-MCs. Takine~\cite{Taki16} presented an algorithm for the conditional stationary distribution vector $\vc{\pi}^{(n)}$ under some additional conditions, which are removed by Kimura and Takine~\cite{M.Kimu18}. As mentioned in Section~1 of \cite{Taki16}, the conditional stationary distribution vector $\vc{\pi}^{(n)}$ can be a good approximation to the stationary distribution vector $\vc{\pi}$ for all sufficiently large $n \in \bbN$. Klimenok and Dudin~\cite{Klim06}, Li et al.~\cite{LiQuan05-RG}, and Shin and Pearce~\cite{Shin98} proposed respective algorithms for approximately computing the stationary distribution vector $\vc{\pi}$ in UBH-MCs by making transition rates (or transition probabilities) eventually {\it level independent}. 

In summary, all the existing algorithms mentioned above are originally designed for approximately computing the stationary distribution vector $\vc{\pi}$ in UBH-MCs (including LD-QBDs). Thus, although those algorithms fulfill {\it practical} purposes, they do not make any {\it theoretical} contribution to the quasi-algorithmic construction of the stationary distribution vector $\vc{\pi}$.

Unlike those previous studies, Masuyama~\cite{Masu19-QUESTA} has taken the first step to the quasi-algorithmic construction of $\vc{\pi}$ in UBH-MCs. To describe the numerical procedure proposed in \cite{Masu19-QUESTA} (and to facilitate later discussion), we introduce the notation: A finite dimensional matrix (resp. vector) is extended (if necessary) to an infinite dimensional matrix (resp. vector) by appending zeros to it with its original elements fixed  in their original positions. This notation enables us to define {\it addition}, {\it subtraction}, and {\it product} over vectors and matrices with different sizes. Furthermore, for $k,n \in \{0,\pm1, \pm2,\dots\}$, the notation ``$\prod_{m=k}^n \downarrow$" denotes a product operator of matrices (including scalars) such that 
\[
\prod_{m=k}^n \downarrow \vc{A}_m 
=
\left\{
\begin{array}{ll}
\vc{A}_n\vc{A}_{n-1} \cdots \vc{A}_k, & k \le n,
\\
\vc{I}, & k >n.
\end{array}
\right.
\]

Masuyama~\cite{Masu19-QUESTA} proposed a {\it sequential update algorithm} (Algorithm~1 therein) that generates a {\it matrix-infinite-product form (MIP-form) solution} for $\vc{\pi}=(\vc{\pi}_0,\vc{\pi}_1,\dots)$:
\begin{align}
\vc{\pi}
=
\lim_{n\to\infty}
\left(
{
\vc{\alpha}_n^{\dag} \vc{U}_n^*
\prod_{m=k}^{n-1} \downarrow \vc{U}_m
\over
\vc{\alpha}_n^{\dag} \vc{U}_n^*
\sum_{\ell=0}^n
(\prod_{m=\ell}^{n-1} \downarrow \vc{U}_m) \vc{e}
}
\right)_{k\in\bbZ_{[0,n]}},
\label{Masu19-MIP-whole}
\end{align}
where (i) $\vc{U}_n = \vc{Q}_{n+1,n}\vc{U}_n^*$; (ii) each $\vc{U}_n^*$, $n\in\bbZ_+$ is computed in a finite number of operations starting from $\vc{U}_0^*=(-\vc{Q}_{0,0})^{-1}$; and (iii) $\vc{\alpha}^{\dag}_n$ is a $1 \times M_n$ probability vector that is an optimal solution for a certain linear fractional programing (LFP) problem, Problem~\ref{prob-01} in Section~\ref{subsec:existing_MIP-form_solution}.  For the reader's convenience, we summarize the sequential update algorithm \cite[Algorithm~1]{Masu19-QUESTA} in Algorithm~\ref{algo-01} in Section~\ref{sec-discussion}.

The existing MIP-form solution (\ref{Masu19-MIP-whole}) (and thus the sequential update algorithm in \cite{Masu19-QUESTA}) requires Conditions~\ref{cond-01} and \ref{cond-02} below (see \cite[Conditions 1 and 2]{Masu19-QUESTA}). 
\begin{CON}[Foster-Lyapunov drift condition]\label{cond-01}
There exist a constant $b \in (0,\infty)$, a finite set $\bbC \subset \bbS$, and a positive column vector $\vc{v}:=(v(k,i))_{(k,i)\in\bbS}$ such that $\inf_{(k,i)\in\bbS}v(k,i) > 0$ and
\begin{equation}
\vc{Q}\vc{v} \le  - \vc{e} + b \vc{1}_{\bbC},
\label{ineqn-QV<=-f+b1_C}
\end{equation}
where
$\vc{1}_{\bbA}:=(1_{\bbA}(k,i))_{(k,i)\in\bbS}$, $\bbA \subseteq \bbS$, denotes a column vector given by
\[
1_{\bbA}(k,i)
=\left\{
\begin{array}{l@{~~~}l}
1, & (k,i) \in \bbA,
\\
0, & (k,i) \in \bbS\setminus\bbA.
\end{array}
\right.
\]
The space of parameter sets $(\vc{v},b,\bbC)$ satisfying (\ref{ineqn-QV<=-f+b1_C}) is denoted by $\cal{V}$, for convenience. 
\end{CON}

\begin{CON}[Convergence condition for the sequential update algorithm in \cite{Masu19-QUESTA}]\label{cond-02}
\begin{eqnarray*}
\vc{\pi}|\vc{Q}| \vc{v}
=
\sum_{(n,i) \in \bbS} \pi(n,i) |q(n,i;n,i)| v(n,i) < \infty,
\end{eqnarray*}
where, for any matrix $\vc{X}$ (resp.~vector $\vc{x}$), the symbol $|\vc{X}|$ (resp.~$|\vc{x}|$) denotes the matrix (resp.~vector) obtained by taking the absolute value of each element of $\vc{X}$ (resp. $\vc{x}$).
\end{CON}

Conditions~\ref{cond-01} and \ref{cond-02} diminish the usefulness of the MIP-form solution (\ref{Masu19-MIP-whole}). A parameter set $(\vc{v},b,\bbC) \in \cal{V}$ of Condition~\ref{cond-01} is required to describe the objective function of Problem~\ref{prob-01} with optimal solution $\vc{\alpha}^{\dag}_n$. In fact, Condition~\ref{cond-01} holds if and only if the irreducible generator $\vc{Q}$ is ergodic (see \cite[Theorem 1.1]{Kont16}). However, even though we know that $\vc{Q}$ is ergodic, it may not be easy to find $(\vc{v},b,\bbC) \in \cal{V}$ in some cases. In addition, Condition~\ref{cond-02} is required to prove the convergence of the MIP-form solution (\ref{Masu19-MIP-whole}) (see \cite[Theorem~3.2]{Masu19-QUESTA}), though this condition does not necessarily hold for all UBH-MCs (such an example is presented in Appendix~\ref{append-cond-02}). 

We note that the existing MIP-form solution (\ref{Masu19-MIP-whole}) is not quasi-algorithmically constructible. For each $n \in \bbZ_+$, the MIP-form solution (\ref{Masu19-MIP-whole}) requires an optimal solution $\vc{\alpha}_n^{\dag}$ for Problem~\ref{prob-01}, and its objective function includes the {\it infinite sums}:
\begin{equation}
\sum_{\ell=n+1}^{\infty} \vc{Q}_{k,\ell}\vc{v}_{\ell},
\qquad n\in\bbZ_+,~k \in \bbZ_{[0,n]},
\label{eqn-infinite-sums}
\end{equation}
where $\vc{v}_{\ell} =(v(\ell,i))_{i\in\bbM_{\ell}}$ is a subvector of $\vc{v}$. Clearly, the infinite sum has {\it infinite} computational complexity, in general. Therefore, the existing MIP-form solution (\ref{Masu19-MIP-whole}) is not quasi-algorithmically constructible, and computing it requires the truncation of the infinite sums (\ref{eqn-infinite-sums}), which causes {\it truncation error}.

This paper presents a new MIP-form solution for $\vc{\pi}=(\vc{\pi}_0,\vc{\pi}_1,\dots)$ quasi-algorithmically constructible:
\begin{equation}
\vc{\pi}
:= \lim_{\scriptstyle n\to\infty \atop \scriptstyle n\ge K}
(\presub{(n)}\wh{\vc{\pi}}^*)_{k \in \bbZ_{[0,n]}}
:= \lim_{\scriptstyle n\to\infty \atop \scriptstyle n\ge K}
\left(
{
\vc{\alpha}_n^* \vc{U}_n^* \prod_{m=k}^{n-1} \downarrow \vc{U}_m
\over
\vc{\alpha}_n^* \vc{U}_n^*
\sum_{\ell=0}^n (\prod_{m=\ell}^{n-1} \downarrow \vc{U}_m) \vc{e}
}
\right)_{\!\!k\in\bbZ_{[0,n]}},
\label{defn-new-MIP-whole}
\end{equation}
where $K \in \bbZ_+$ is a free parameter, and where $\vc{\alpha}_n^*$, $n\ge K$,  is a $1 \times M_n$ probability vector that is an optimal solution for the LFP problem (Problem~\ref{prob-02} in Section~\ref{sec-main}) different from Problem~\ref{prob-01} of the existing MIP-form solution (\ref{Masu19-MIP-whole}). The objective function of Problem~\ref{prob-02} does not include such a parameter set as $(\vc{v},b,\bbC) \in \cal{V}$ and is computed in a finite number of operations without any infinite sum like (\ref{eqn-infinite-sums}). Thus, the new MIP-form solution (\ref{defn-new-MIP-whole}) is quasi-algorithmically constructible. In addition, the new solution holds without any additional conditions, such as Condition~\ref{cond-02} for convergence, except for the ergodicity of~$\vc{Q}$.


The rest of this paper is divided into five sections. Section~\ref{sec-preliminary} provides preliminary results including the existing MIP-form solution presented in \cite{Masu19-QUESTA}. Section~\ref{sec-main} presents our new MIP-form solution together with its theoretical foundations, which is the main result of this paper. Section~\ref{sec-discussion} considers the advantages of our MIP-form solution over the existing one. Section~\ref{sec-concluding} discusses special cases where our solution can be established more effectively. Finally, Section~\ref{sec:concluding-remarks} contains concluding remarks.


\section{The existing MIP-form solution}\label{sec-preliminary}

This section reviews the existing MIP-form solution (\ref{Masu19-MIP-whole}) to facilitate describing our new MIP-form solution (established in Section~\ref{sec-main}) and understanding its advantages over the existing one. First, we introduce the foundation of the existing MIP-form solution, the {\it last-block-column linearly augmented (LBC-LA)} truncation approximation to the stationary distribution vector $\vc{\pi}$. We then provide a matrix-product form of the LBC-LA truncation approximation. Finally, as a certain limit of the matrix-product form, we express the existing MIP-form solution (\ref{Masu19-MIP-whole}) for $\vc{\pi}$.

For later reference, we define the $\ell_1$-distance for vectors of different dimensions. Let $\bbG$ be a countable (possibly infinite) set. Let $\vc{h}_1 =(h_1(j))_{j\in\bbG_1}$ and $\vc{h}_2:=(h_2(j))_{j\in\bbG_2}$ be real-valued vectors, where $\bbG_1$ and $\bbG_2$ are subsets of $\bbG$. For vectors $\vc{h}_1$ and $\vc{h}_2$, let
\[
\| \vc{h}_1 - \vc{h}_2 \|_1
:= 
\sum_{j \in \bbG_1 \cap \bbG_2} | h_1(j) - h_2(j) |
+ \sum_{j \in \bbG_1\setminus\bbG_2} | h_1(j) |
+ \sum_{j \in \bbG_2\setminus\bbG_1} | h_2(j) |,
\]
which is referred to as the $\ell_1$-distance of the difference between $\vc{h}_1$ and $\vc{h}_2$.

\subsection{Definition of the LBC-LA truncation approximation}

We begin with some definitions to describe the LBC-LA truncation approximation. Let $\bbS_n = \bigcup_{k=0}^n \bbL_k$ for $n \in \bbZ_+$ and $\bbS_{-1} = \varnothing$ for convenience. For any set $\bbX$, let $|\bbX|$ denote the cardinality of $\bbX$. For any $n \in \bbZ_+$, let $\presub{(n)}\vc{Q}$ denote a submatrix of $\vc{Q}$ consisting of the first $|\bbS_n|$ rows and columns, that is,
%
%
\begin{eqnarray*}
\presub{(n)}\vc{Q}
=
\bordermatrix{
		&
\bbL_0 	&
\bbL_1 	&
\bbL_2 	&
\cdots		&
\bbL_{n-2} 	&
\bbL_{n-1} 	&
\bbL_n 		
\cr
\bbL_0		&
\vc{Q}_{0,0} 	&
\vc{Q}_{0,1} 	&
\vc{Q}_{0,2} 	&
\cdots		&
\vc{Q}_{0,n-2} 	&
\vc{Q}_{0,n-1} 	&
\vc{Q}_{0,n} 	
\cr
\bbL_1		&
\vc{Q}_{1,0} 	&
\vc{Q}_{1,1} 	&
\vc{Q}_{1,2} 	&
\cdots		&
\vc{Q}_{1,n-2} 	&
\vc{Q}_{1,n-1} 	&
\vc{Q}_{1,n} 	
\cr
\bbL_2		&
\vc{O} 		&
\vc{Q}_{2,1} 	&
\vc{Q}_{2,2} 	&
\cdots		&
\vc{Q}_{2,n-2} 	&
\vc{Q}_{2,n-1} 	&
\vc{Q}_{2,n} 	
\cr
~\vdots		&
\vdots		&
\vdots		&
\vdots		&
\ddots		&
\vdots		&
\vdots		&
\vdots		
\cr
\bbL_{n-1}	&
\vc{O}		&
\vc{O}		&
\vc{O}		&
\cdots		&
\vc{Q}_{n-1,n-2} 	&
\vc{Q}_{n-1,n-1} 	&
\vc{Q}_{n-1,n} 
\cr
\bbL_n		&
\vc{O}		&
\vc{O}		&
\vc{O}		&
\cdots		&
\vc{O}		&
\vc{Q}_{n,n-1} 	&
\vc{Q}_{n,n} 
}.
\end{eqnarray*}
%
Furthermore, let $\presub{(n)}\wh{\vc{\alpha}}$, $n \in \bbZ_+$, denote a $1 \times |\bbS_n|$ probability vector that has its probability masses on the last block (corresponding to $\bbL_n$), that is, 
\begin{equation*}
\presub{(n)}\wh{\vc{\alpha}}
=
\bordermatrix{
        	& 
\bbS_{n-1} &
\bbL_n 
\cr
			& 
\vc{0}  	& 
\vc{\alpha}_n	
},
\end{equation*}
where $\vc{\alpha}_n:=(\alpha_n(j))_{j\in\bbM_n}$ is a probability vector.

We define the LBC-LA truncation approximation to $\vc{\pi}$ as a stationary distribution vector of the finite essential $Q$-matrix $\presub{(n)}\wh{\vc{Q}}$ given by
\begin{equation*}
\presub{(n)}\wh{\vc{Q}}
= \presub{(n)}\vc{Q} 
- \presub{(n)}\vc{Q}\vc{e} \presub{(n)}\wh{\vc{\alpha}},
\qquad n\in\bbZ_+.
\end{equation*}
The $Q$-matrix $\presub{(n)}\wh{\vc{Q}}$ is referred to as the {\it last-block-column linearly augmented (LBC-LA) truncation} of generator $\vc{Q}$. Clearly, the
LBC-LA truncation $\presub{(n)}\wh{\vc{Q}}$ has at least one stationary distribution vector. One of them is
given by
\begin{equation*}
\presub{(n)}\wh{\vc{\pi}} 
= {\presub{(n)}\wh{\vc{\alpha}} (- \presub{(n)}\vc{Q})^{-1}
\over \presub{(n)}\wh{\vc{\alpha}} (- \presub{(n)}\vc{Q})^{-1}\vc{e}},
\qquad n \in \bbZ_+.
\end{equation*}
This stationary distribution vector $\presub{(n)}\wh{\vc{\pi}}$ is referred to as the LBC-LA truncation approximation to $\vc{\pi}$, and is partitioned level-wise:
\begin{align*}
\presub{(n)}\wh{\vc{\pi}}
=
\bordermatrix{
       	& 
\bbL_0 	& 
\bbL_1 	& 
\cdots	&
\bbL_n 	
\cr
 	& 
\presub{(n)}\wh{\vc{\pi}}_0 & 
\presub{(n)}\wh{\vc{\pi}}_1 & 
\cdots				&	
\presub{(n)}\wh{\vc{\pi}}_n
}.
\end{align*}

\subsection{A matrix-product form of the LBC-LA truncation approximation}

This subsection presents a matrix-product form of the LBC-LA truncation approximation $\presub{(n)}\wh{\vc{\pi}}$. We first introduce some vectors and matrices needed to express the matrix-product form of $\presub{(n)}\wh{\vc{\pi}}$. We then present the matrix-product form together with recursive formulas for the components of it. Finally, we mention a connection of the matrix-product form to the existing MIP-form solution (\ref{Masu19-MIP-whole}) for $\vc{\pi}$.

We define the component matrices of the matrix-product form of $\presub{(n)}\wh{\vc{\pi}}$. Let $\vc{U}_n^*$, $n \in \bbZ_+$, denote an $M_n \times M_n$ matrix such that its $(i,j)$-th element represents the expected total sojourn time in state $(n,j)$ before the first visit to $\ol{\bbS}_n:=\bbS\setminus \bbS_n=\cup_{k=n+1}^{\infty}\bbL_k$ (i.e., to any state in levels $n+1,n+2,\dots$) starting from state $(n,i)$. The matrices $\vc{U}_n^*$, $n \in \bbZ_+$ are determined recursively:
\begin{align}
\vc{U}_0^*
&= (-\vc{Q}_{0,0})^{-1},
\label{recursion:U_n^*-01}
\end{align}
and, for $n \in \bbN$,
\begin{align}
\vc{U}_n^*
&= \left[
-\vc{Q}_{n,n} 
- 
\dm\sum_{\ell=0}^{n-1}
(\vc{Q}_{n,n-1}\vc{U}_{n-1}^*)
(\vc{Q}_{n-1,n-2}\vc{U}_{n-2}^*) \cdots
(\vc{Q}_{\ell+1,\ell}\vc{U}_{\ell}^*) 
\vc{Q}_{\ell,n}
\right]^{-1}.
\label{recursion:U_n^*-02}
\end{align}
Furthermore, let $\vc{U}^*_{n,k}$, $n \in \bbZ_+$, $k\in\bbZ_{[0,n]}$ denote
\begin{eqnarray}
\vc{U}^*_{n,k}
&=&  \vc{U}_n^* 
\cdot
(\vc{Q}_{n,n-1}\vc{U}_{n-1}^*)
(\vc{Q}_{n-1,n-2}\vc{U}_{n-2}^*) \cdots
(\vc{Q}_{k+1,k}\vc{U}_k^*).
\label{defn-U_{n,k}^*}
\end{eqnarray}
It then follows from (\ref{recursion:U_n^*-02}) that
\begin{eqnarray}
\vc{U}_n^*
&=& 
\left(
- \vc{Q}_{n,n} 
- \vc{Q}_{n,n-1} \sum_{\ell=0}^{n-1} \vc{U}_{n-1,\ell}^*\vc{Q}_{\ell,n}
\right)^{-1},
\qquad n \in \bbN.
\label{eqn-U_n^*}
\end{eqnarray}
For simplicity of notation, we define $\vc{U}_m$, $m\in\bbZ_+$ as
\begin{align}
\vc{U}_m = \vc{Q}_{m+1,m}\vc{U}_m^*.
\label{defn:U_m}
\end{align}
Using (\ref{defn:U_m}), we rewrite (\ref{defn-U_{n,k}^*}) as 
\begin{eqnarray}
\vc{U}^*_{n,k}
&=&  \vc{U}_n^* \prod_{m=k}^{n-1} \downarrow \vc{U}_m,
\qquad k\in\bbZ_{[0,n]},
\label{eqn:U_{n,k}^*}
\end{eqnarray}
Finally, $\vc{u}^*_n:=(u_n^*(i))_{i\in\bbM_n}$, $n \in \bbZ_+$, is defined as the row sum vector of the $|\bbL_n| \times |\bbS_n|$ matrix $(\vc{U}_{n,0}^*~\vc{U}_{n,1}^*~\cdots~\vc{U}_{n,n}^*)$:
\begin{eqnarray}
\vc{u}_n^* 
&=& \sum_{\ell=0}^n \vc{U}_{n,\ell}^*\vc{e} 
= \sum_{\ell=0}^n \vc{U}_n^*
\left(\prod_{m=\ell}^{n-1} \downarrow \vc{U}_m \right) \vc{e}
>\vc{0},
\label{defn-u_n^*}
\end{eqnarray}
where $\vc{u}_n^* > \vc{0}$ follows from \cite[Remark~2.3]{Masu19-QUESTA}. 

\begin{REM}\label{rem-U_{n,k}^*}
The probabilistic interpretation of $\vc{U}_n^*$ leads to that of $\vc{U}_{n,\ell}^*$: Its $(i,j)$-th element represents the expected total sojourn time in state $(\ell,j)$ before the first visit to $\ol{\bbS}_{n}$ starting from state $(n,i)$. This interpretation implies that the matrix $(\vc{U}_{n,0}^*~\vc{U}_{n,1}^*~\cdots~\vc{U}_{n,n}^*)$ is identified with the rows in level $n$ of $(-\presub{(n)}\vc{Q})^{-1}$ because this inverse matrix consists of the expected total sojourn times in the respective states in $\bbS_n$ before the first visit to $\ol{\bbS}_{n}$ starting from the states in $\bbS_n$. We thus have (see \cite[Theorem~5.5.1]{Lato99} and \cite[Remark~2.2]{Masu19-QUESTA}) 
\begin{align}
\vc{\pi}_{\ell} 
&= \vc{\pi}_{n+1} \vc{Q}_{n+1,n}\vc{U}_{n,\ell}^*
\nonumber
\\
&=\vc{\pi}_{n+1} \vc{U}_n\vc{U}_{n-1} \cdots \vc{U}_{\ell},
\qquad n \in \bbZ_+,~\ell \in \bbZ_{[0,n]},
\label{eqn:210618-01}
\end{align}
where the second equality is due to (\ref{defn:U_m}) and (\ref{eqn:U_{n,k}^*}). \end{REM}

We are now ready to present the matrix-product form of $\presub{(n)}\wh{\vc{\pi}}=(\presub{(n)}\wh{\vc{\pi}}_0,\presub{(n)}\wh{\vc{\pi}}_1,\dots,\presub{(n)}\wh{\vc{\pi}}_n)$.
\begin{PROP}[{}{\cite[Equation~(2.21) and Lemma~2.2]{Masu19-QUESTA}}]\label{lem-product-form-(s)ol{pi}_s}
For $n \in \bbZ_+$ and $k \in \bbZ_{[0,n]}$,
\begin{eqnarray}
\presub{(n)}\wh{\vc{\pi}}_{k}
&=&
{
\vc{\alpha}_n \vc{U}_{n,k}^*
\over 
\vc{\alpha}_n \vc{u}_n^*
}
=
{
\vc{\alpha}_n \vc{U}_n^* \prod_{m=k}^{n-1} \downarrow \vc{U}_m
\over 
\vc{\alpha}_n \vc{U}_n^*
\sum_{\ell=0}^n (\prod_{m=\ell}^{n-1} \downarrow \vc{U}_m) \vc{e}
}.
\label{eqn-(n_s)pi_{n_s,l}}
\end{eqnarray}
\end{PROP}

The components $\vc{U}_{n,k}^*$ and $\vc{u}_n^*$ of the matrix-product form (\ref{eqn-(n_s)pi_{n_s,l}}) are the following recursive formulas  (see \cite[Equations~(3.23) and (3.24)]{Masu19-QUESTA}):
\begin{subequations}\label{recursion-U_{n,k}^*}
\begin{eqnarray}
\vc{U}_{0,0}^*
&=& \vc{U}_0^* = (-\vc{Q}_{0,0})^{-1},
\label{eqn-U_{0,0}^*}
\\
\vc{U}_{n,k}^*
&=& 
\left\{
\begin{array}{l@{~~~}ll}
\vc{U}_n^*,
& n \in \bbN,\ & k = n,
\\
\vc{U}_n^* \vc{Q}_{n,n-1} \cdot \vc{U}^*_{n-1,k},
& n \in \bbN,\ & k=n-1,n-2,\dots,0,
\\
\end{array}
\right.
\label{eqn-U_{n,k}^*}
\end{eqnarray}
\end{subequations}
and
\begin{subequations}\label{recursion-u_n^*}
\begin{eqnarray}
\vc{u}^*_0
&=& \vc{U}_0^* \vc{e} = (-\vc{Q}_{0,0})^{-1}\vc{e},
\label{eqn-u_0^*}
\\
\vc{u}^*_n
&=& \vc{U}_n^* 
\left( \vc{e} + \vc{Q}_{n,n-1} \vc{u}^*_{n-1}
\right),
\qquad n \in \bbN,
\label{eqn-u_n^*}
\end{eqnarray}
\end{subequations}
where $\vc{U}_n^*$ is given in (\ref{eqn-U_n^*}).

Proposition~\ref{lem-product-form-(s)ol{pi}_s} implies that an MIP-form solution for $\vc{\pi}$ can be obtained as a limit of the matrix-product form (\ref{eqn-(n_s)pi_{n_s,l}}), that is,
\begin{align*}
\left( 
{
\vc{\alpha}_n \vc{U}_n^* \prod_{m=k}^{n-1} \downarrow \vc{U}_m
\over 
\vc{\alpha}_n \vc{U}_n^*
\sum_{\ell=0}^n (\prod_{m=\ell}^{n-1} \downarrow \vc{U}_m) \vc{e}
}
\right)_{k \in \bbZ_{[0,n]}} \to \vc{\pi}\quad \mbox{as $n \to \infty$}.
\end{align*}
However, as mentioned in \cite[Section~2.3]{Masu19-QUESTA}, achieving such a limit convergent to $\vc{\pi}$ requires choosing the probability vectors $\vc{\alpha}_n$, $n \in \bbZ_+$ appropriately. To do this, Masuyama~\cite{Masu19-QUESTA} formulated a series of certain linear fractional programing (LFP) problems. The details are given in the next subsection.

\subsection{The existing MIP-form solution for the stationary distribution vector}\label{subsec:existing_MIP-form_solution}

This subsection provides a brief summary of the existing MIP-form solution (\ref{Masu19-MIP-whole}) for $\vc{\pi}$. We first show the LFP problem whose optimal solution is the key component of the existing MIP-form solution. We then present the proposition that summarizes the theoretical results on the existing MIP-form solution.

To obtain an MIP-form solution for $\vc{\pi}$, Masuyama~\cite{Masu19-QUESTA} formulates the following LFP problem indexed with $n \in \bbZ_+$:
\begin{PROB}\label{prob-01}
\begin{subequations}
\begin{alignat*}{2}
& \mbox{\rm Minimize} 
&
& {\vc{\alpha}_n\vc{y}_n  \over \vc{\alpha}_n\vc{u}_n^*};
\\
& \mbox{\rm Subject to}	\quad  
&  	
& 
\vc{\alpha}_n\vc{e} = 1,
\\	
&		
&
&
\vc{\alpha}_n \ge \vc{0},
\end{alignat*}
\end{subequations}
where $\vc{y}_n:=(y_n(i))_{i\in\bbM_n}$ is given by
\begin{equation}
\vc{y}_n = 
\vc{v}_n 
+ 
\sum_{k=0}^n \vc{U}_{n,k}^*
\sum_{\ell=n+1}^{\infty} \vc{Q}_{k,\ell} \vc{v}_{\ell} > \vc{0}.
\label{defn-y_n}
\end{equation}
\end{PROB}

\begin{REM}\label{rem:objective_function}
Equation (\ref{defn-u_n^*}) implies that $\vc{\alpha}_n\vc{u}_n^* > 0$ for any feasible solution $\vc{\alpha}_n$. Therefore, the objective function of Problem~\ref{prob-01} is well-defined. 
\end{REM}

\begin{REM}\label{rem-LD-QBD}
If the UBH generator $\vc{Q}$ satisfies $\vc{Q}_{k,k+m} = \vc{O}$ for $m \ge 2$ ($\vc{Q}$ is an LD-QBD generator), then (\ref{defn-y_n}) is reduced to
\[
\vc{y}_n = \vc{v}_n + \vc{U}_{n,n}^*\vc{Q}_{n,n+1} \vc{v}_{n+1}.
\]
\end{REM}

Proposition~\ref{prop-(n)wh{pi}^{dag}} below shows how to find an optimal solution for Problem~\ref{prob-01} and how to construct, with such a solution, the existing MIP-form solution (\ref{Masu19-MIP-whole}), originally presented in \cite{Masu19-QUESTA}.
\begin{PROP}[{}{\cite[Theorems~3.1 and 3.2, and Corollary~3.1]{Masu19-QUESTA}}]\label{prop-(n)wh{pi}^{dag}}
If Condition~\ref{cond-01} is satisfied, then the following hold.
\begin{enumerate}
\item For $n \in \bbZ_+$, let 
 $\vc{\alpha}_n^{\dag}:=(\alpha_n^{\dag}(j))_{j\in\bbM_n}$ denote a probability vector such that
\begin{equation*}
\alpha_n^{\dag}(j)
=
\left\{
\begin{array}{ll}
1, & \quad j=j_n^{\dag},
\\
0, & \quad j \neq j_n^{\dag},
\end{array}
\right.
\end{equation*}
where
\begin{equation}
j_n^{\dag} 
\in \argmin_{j\in\bbM_n} {y_n(j)  \over u_n^*(j)}.
\label{defn-j_n^{dag}}
\end{equation}
The probability vector $\vc{\alpha}_n^{\dag}$ is an optimal solution for Problem~\ref{prob-01}.
\item For $n \in \bbZ_+$,  
let 
$\presub{(n)}\wh{\vc{\pi}}^{\dag}
:= (
\presub{(n)}\wh{\vc{\pi}}_{0}^{\dag},
\presub{(n)}\wh{\vc{\pi}}_{1}^{\dag},
\dots,
\presub{(n)}\wh{\vc{\pi}}_{n}^{\dag})$ denote a probability vector such that
\begin{equation}
\presub{(n)}\wh{\vc{\pi}}_{k}^{\dag}
= {
\vc{\alpha}_n^{\dag}\vc{U}_{n,k}^*
\over
\vc{\alpha}_n^{\dag} \vc{u}_n^*
}
= {
\row_{j_n^{\dag}}(\vc{U}_{n,k}^*)
\over  
u_n^*(j_n^{\dag}) 
},\qquad k \in \bbZ_{[0,n]},
\label{defn-(n)wh{pi}_k^{dag}}
\end{equation}
where $\row_j(\,\cdot\,)$ denotes the $j$-th row of the matrix in the parentheses. Furthermore, if a parameter set $(\vc{v},b,\bbC) \in \cal{V}$ of Condition~\ref{cond-01} satisfies Condition~\ref{cond-02}, then
\begin{eqnarray*}
\lim_{n\to\infty} \| \presub{(n)}\wh{\vc{\pi}}^{\dag} - \vc{\pi} \|_1  &=& 0,
\end{eqnarray*}
and
\begin{equation}
\vc{\pi}
=\lim_{n\to\infty}
\left(
{ 
\vc{\alpha}_n^{\dag} \vc{U}_{n,k}^*
\over 
\vc{\alpha}_n^{\dag} \vc{u}_n^*
}
\right)_{\!\!k\in\bbZ_{[0,n]}}
= \lim_{n\to\infty}
\left(
{ 
\vc{\alpha}_n^{\dag} \vc{U}_n^* \prod_{m=k}^{n-1} \downarrow \vc{U}_m
\over 
\vc{\alpha}_n^{\dag} \vc{U}_n^*
\sum_{\ell=0}^n (\prod_{m=\ell}^{n-1} \downarrow \vc{U}_m) \vc{e}
}
\right)_{\!\!k\in\bbZ_{[0,n]}}.
\label{eqn-MIP-Masu19-QUESTA}
\end{equation}
\end{enumerate}

\end{PROP}

\begin{REM}
The symbols $\vc{\alpha}_n^{\dag}$, $j_n^{\dag}$, and $\presub{(n)}\wh{\vc{\pi}}^{\dag}$ correspond to $\vc{\alpha}_n^*$\, $j_n^*$, and $\presub{(n)}\wh{\vc{\pi}}^*$, respectively, in the original notation of \cite{Masu19-QUESTA}.
\end{REM}

\begin{REM}\label{rem:old-MIP-form}
Constructing the MIP-form solution (\ref{eqn-MIP-Masu19-QUESTA}) (and thus (\ref{Masu19-MIP-whole})) requires the sequence $\{\vc{\alpha}_n^{\dag}\}$. Each $\vc{\alpha}_n^{\dag}$ is an optimal solution for Problem~\ref{prob-01} with index $n$. Therefore, we only have to solve a single  optimization problem for each $n \in \bbZ_+$.
\end{REM}


\section{A new MIP-form solution}\label{sec-main}

This section presents a new MIP-form solution for $\vc{\pi}$. To construct the new MIP-form solution, we first formulate an LFP problem different from Problem~\ref{prob-01}. We then provide a way to find its optimal solution. Finally, we show the new MIP-form solution for $\vc{\pi}$.

We introduce some definitions to describe our LFP problem for the new MIP-form solution. Let $\bbZ_{\ge m}:=\{m,m+1,\dots\}$ for any integer $m$. We then define $\vc{u}_{n,\bbK}^*:=(u_{n,\bbK}^*(i))_{i\in\bbM_n}$, $n \in \bbZ_{\ge K}$, as
\begin{equation}
\vc{u}_{n,\bbK}^* 
= \sum_{\ell \in \bbK} \vc{U}_{n,\ell}^*\vc{e} 
= \vc{U}_n^* \sum_{\ell \in \bbK} 
\left(
\prod_{m=\ell}^{n-1}\downarrow \vc{U}_m
\right) \vc{e},
\label{defn-u_{n,N}^*}
\end{equation}
where the second equality is due to (\ref{eqn:U_{n,k}^*}). We also define $\bbI_n^+$, $n \in \bbZ_+$, as
\begin{eqnarray}
\bbI_n^+ 
&=& 
\left\{
j\in\bbM_n: (\vc{e}^{\top}\vc{Q}_{n+1,n})_j > 0
\right\},\qquad n \in \bbZ_+,
\label{defn:I_n^+}
\end{eqnarray}
where $(\,\cdot\,)_j$ denotes the $j$-th element of the (row or column) vector in the parentheses. Since the generator $\vc{Q}$ is irreducible, $(\vc{e}^{\top}\vc{Q}_{n+1,n})_j > 0$ for at least one $j \in \bbM_n$ and thus $\bbI_n^+ \neq \varnothing$. Furthermore, since $\vc{\pi}_{n+1} > \vc{0}$ and $\vc{Q}_{n+1,n} \ge \vc{O}, \neq \vc{O}$,
\begin{align}
\left\{
\begin{array}{ll}
(\vc{\pi}_{n+1}\vc{Q}_{n+1,n})_j > 0, & j \in \bbI_n^+,
\\
(\vc{\pi}_{n+1}\vc{Q}_{n+1,n})_j = 0, & j \in \bbM_n\setminus \bbI_n^+,
\end{array}
\right.
\label{eqn-e*Q_{n+1,n}}
\end{align}
which is used later.

The following is our LFP problem, which is formulated for each $n \in \bbZ_{\ge K}$.
\begin{PROB}\label{prob-02}
%
\begin{alignat}{2}
& \mbox{\rm Maximize} 
&
\quad r_n(\vc{\alpha}_n)&:={\vc{\alpha}_n\vc{u}_{n,\bbK}^*  \over \vc{\alpha}_n\vc{u}_n^*};
\label{defn-r_n} 
\\
& \mbox{\rm Subject to}	
&  	
\vc{\alpha}_n\vc{e} \!&\phantom{:}= 1,
\label{prob-02-const-xe=1}
\\
&		
&
\alpha_n(j)\!&\phantom{:}=0,\quad j \in \bbM_n\setminus\bbI_n^+,
\label{prob-02-const-alpha(j)=0}
\\	
&		
&
\vc{\alpha}_n \!&\phantom{:}\ge \vc{0}.
\label{prob-02-const-x>=0}
\end{alignat}
%
%
\end{PROB}

\begin{REM}
The objective function of Problem~\ref{prob-02} is well-defined, as with Problem~\ref{prob-01} (see Remark~\ref{rem:objective_function}).
\end{REM}

An optimal solution for Problem~\ref{prob-02} is given by Lemma~\ref{lem-optima-solution-LFP} below.
\begin{LEM}\label{lem-optima-solution-LFP}
For $n \in \bbZ_{\ge K}$, fix $j_n^* \in \bbI_n^+ \neq \varnothing$ such that
\begin{equation}
j_n^* \in \bbJ_n^* 
:= \argmax_{j\in\bbI_n^+} {u_{n,\bbK}^*(j)  \over u_n^*(j)}
\subseteq \bbI_n^+,
\label{defn-J_n^*}
\end{equation}
and let $\vc{\alpha}_n^*:=(\alpha_n^*(j))_{j\in\bbM_n}$ denote a unit row vector such that
\begin{equation}
\alpha_n^*(j)
=
\left\{
\begin{array}{l@{~~~}l}
1, & j=j_n^*,
\\
0, & j \neq j_n^*.
\end{array}
\right.
\label{defn:alpha_n^*}
\end{equation}
The vector $\vc{\alpha}_n^*$ is then an optimal solution for Problem~\ref{prob-02}, and the optimal (and thus maximum) value $r_n(\vc{\alpha}_n^*)=u_{n,\bbK}^*(j_n^*)/u_{n}^*(j_n^*)$ is positive.
\end{LEM}

\begin{proof}
First, we prove that $\vc{\alpha}_n^*$ is an optimal solution for Problem~\ref{prob-02}. It follows from (\ref{defn-J_n^*}) and (\ref{defn:alpha_n^*}) that
\begin{equation}
\xi_n := {\vc{\alpha}_n^* \vc{u}_{n,\bbK}^* \over \vc{\alpha}_n^* \vc{u}_n^*}
= {u_{n,\bbK}^*(j_n^*) \over u_n^*(j_n^*)}
= \max_{j\in\bbI_n^+}{u_{n,\bbK}^*(j) \over u_n^*(j)},
\label{defn-xi_n}
\end{equation}
which leads to $u_{n,\bbK}^*(j) \le \xi_n u_n^*(j)$ for all $j \in \bbI_n^+$. Thus, for any feasible solution $\vc{\alpha}_n$ satisfying (\ref{prob-02-const-xe=1})--(\ref{prob-02-const-x>=0}), we have
\begin{equation*}
\vc{\alpha}_n\vc{u}_{n,\bbK}^*
= \sum_{j\in\bbI_n^+}\alpha_n(j) u_{n,\bbK}^*(j) 
\le \xi_n \sum_{j\in\bbI_n^+}\alpha_n(j)u_n^*(j)
= \xi_n \vc{\alpha}_n\vc{u}_n^*.
\end{equation*}
Combining this result with (\ref{defn-r_n}) and (\ref{defn-xi_n}) yields
\[
r_n(\vc{\alpha}_n)
= {\vc{\alpha}_n\vc{u}_{n,\bbK}^* \over \vc{\alpha}_n\vc{u}_n^*} \le \xi_n
= {\vc{\alpha}_n^* \vc{u}_{n,\bbK}^* \over \vc{\alpha}_n^* \vc{u}_n^*}
= r_n(\vc{\alpha}_n^*),
\]
which shows that $\vc{\alpha}_n^*$ is an optimal solution for Problem~\ref{prob-02}.

Next, we prove $r_n(\vc{\alpha}_n^*) > 0$ by contradiction. To this end,  suppose that $r_n(\vc{\alpha}_n^*) = 0$, or equivalently, that $\vc{\alpha}_n \vc{u}_{n,\bbK}^* = 0$ for any $\vc{\alpha}_n$ satisfying (\ref{prob-02-const-xe=1})--(\ref{prob-02-const-x>=0}). We then have
\begin{align*}
u_{n,\bbK}^*(i) = 0\quad \mbox{for all $i \in \bbI_n^+$}.
\end{align*}
Therefore, using (\ref{eqn-e*Q_{n+1,n}}), (\ref{defn-u_{n,N}^*}), and (\ref{eqn:210618-01}), we obtain
\begin{align*}
0 
&= \vc{\pi}_{n+1}\vc{Q}_{n+1,n}\vc{u}_{n,\bbK}^*
=\vc{\pi}_{n+1}\vc{Q}_{n+1,n} \sum_{\ell \in \bbK}  \vc{U}_{n,\ell}^*\vc{e}
= \sum_{\ell \in \bbK} \vc{\pi}_{\ell}\vc{e},
\end{align*}
which contradicts $\vc{\pi}>\vc{0}$. The proof has been completed.
\end{proof}

The next theorem shows that the optimal solution $\vc{\alpha}_n^*$ for Problem~\ref{prob-02} leads to an MIP-form solution for $\vc{\pi}$, which is different from the existing one (\ref{eqn-MIP-Masu19-QUESTA}) based on Problem~\ref{prob-01}.
\begin{THM}\label{thm-E^*-01}
For each $n \in \bbZ_{\ge K}$, let $\presub{(n)}\wh{\vc{\pi}}^*
= (
\presub{(n)}\wh{\vc{\pi}}_{0}^*,
\presub{(n)}\wh{\vc{\pi}}_{1}^*,
\dots,
\presub{(n)}\wh{\vc{\pi}}_{n}^*)$, where $\presub{(n)}\wh{\vc{\pi}}_{k}^*$, $k\in\bbZ_{[0,n]}$, is a $1 \times M_k$ vector obtained by replacing $\vc{\alpha}_n$ in (\ref{eqn-(n_s)pi_{n_s,l}}) with the optimal solution $\vc{\alpha}_n^*$ for Problem~\ref{prob-02}, i.e.,
\begin{equation}
\presub{(n)}\wh{\vc{\pi}}_{k}^*
= {
\vc{\alpha}_n^*\vc{U}_{n,k}^*
\over
\vc{\alpha}_n^* \vc{u}_n^*
}
= {
\row_{j_n^*}(\vc{U}_{n,k}^*)
\over  
u_n^*(j_n^*) 
},\qquad k \in \bbZ_{[0,n]},
\label{defn-(n)wh{pi}_k^*}
\end{equation}
where $j_n^* \in \bbJ_n^*$ defined in (\ref{defn-J_n^*}).
We then have
\begin{eqnarray*}
\lim_{n\to\infty} \| \presub{(n)}\wh{\vc{\pi}}^* - \vc{\pi} \|_1  &=& 0,
\end{eqnarray*}
and therefore 
\begin{eqnarray}
\vc{\pi}
&=& \lim_{n\to\infty}
\left(
{ 
\vc{\alpha}_n^* \vc{U}_{n,k}^*
\over 
\vc{\alpha}_n^* \vc{u}_n^*
}
\right)_{\!\!k\in\bbZ_{[0,n]}}
= \lim_{n\to\infty}
\left(
{ 
\vc{\alpha}_n^* \vc{U}_n^* \prod_{m=k}^{n-1} \downarrow \vc{U}_m
\over 
\vc{\alpha}_n^*  \vc{U}_n^*
\sum_{\ell=0}^n (\prod_{m=\ell}^{n-1} \downarrow \vc{U}_m) \vc{e}
}
\right)_{\!\!k\in\bbZ_{[0,n]}}.\qquad
\label{new-MIP-01}
\end{eqnarray}
\end{THM}

\begin{REM}\label{rem:new-MIP-form}
We only have to solve a single optimization problem (Problem~\ref{prob-02}) for each $n \in \bbZ_+$ to construct the new MIP-form solution (\ref{new-MIP-01}), as with the existing one (\ref{eqn-MIP-Masu19-QUESTA}) (see Remark~\ref{rem:old-MIP-form}).
\end{REM}

{\it Proof of Theorem~\ref{thm-E^*-01}.}~~According to Theorem~\ref{thm-appen-03}, it suffices to show that
\begin{equation}
\liminf_{n\to\infty} 
\sum_{\ell\in\bbK} \presub{(n)}\wh{\vc{\pi}}_{\ell}^* \vc{e}
\ge \sum_{\ell\in\bbK} \vc{\pi}_{\ell}\vc{e}.
\label{lim-r_n(psi_n^*)}
\end{equation}
To do this, we define $\wt{\vc{\alpha}}_n:=(\wt{\alpha}_n(j))_{j\in\bbM_n}$ as
\begin{equation}
\wt{\vc{\alpha}}_n 
= {
\vc{\pi}_{n+1}\vc{Q}_{n+1,n}
\over 
\vc{\pi}_{n+1}\vc{Q}_{n+1,n}\vc{e}
},
\qquad n \in \bbZ_+.
\label{defn-wt{alpha}_n}
\end{equation}
Equations (\ref{eqn-e*Q_{n+1,n}}) and (\ref{defn-wt{alpha}_n}) imply that  
\begin{align*}
\wt{\vc{\alpha}}_n\vc{e} &=1,
\\
\wt{\alpha}_n(j) &= 0,\qquad j \in \bbM_n\setminus\bbI_n^+,
\\
\wt{\vc{\alpha}}_n &\ge \vc{0},
\end{align*}
and thus $\wt{\vc{\alpha}}_n$ is a feasible solution for Problem~\ref{prob-02}. Recall here that $\vc{\alpha}_n^*$ maximizes the objective function $r_n(\,\cdot\,)$ (see Lemma~\ref{lem-optima-solution-LFP}). Therefore, it follows from (\ref{defn-r_n}), (\ref{defn-u_{n,N}^*}), and (\ref{defn-(n)wh{pi}_k^*}) that 
\begin{align}
{\wt{\vc{\alpha}}_n \vc{u}_{n,\bbK}^* \over \wt{\vc{\alpha}}_n\vc{u}_n^*}
=
r_n(\wt{\vc{\alpha}}_n) 
\le r_n(\vc{\alpha}_n^*)
= {\vc{\alpha}_n^*\vc{u}_{n,\bbK}^* \over \vc{\alpha}_n^*\vc{u}_n^*}
= \sum_{\ell\in\bbK} \presub{(n)}\wh{\vc{\pi}}_{\ell}^* \vc{e}.
\label{eqn:210712-01}
\end{align}
It also follows from (\ref{defn-wt{alpha}_n}), (\ref{defn-u_{n,N}^*}), and (\ref{defn-u_n^*}) that
\begin{align}
{\wt{\vc{\alpha}}_n \vc{u}_{n,\bbK}^* \over \wt{\vc{\alpha}}_n\vc{u}_n^*}
&= 
{
\vc{\pi}_{n+1}\vc{Q}_{n+1,n} 
\sum_{\ell\in\bbK} \vc{U}_{n,\ell}^*\vc{e}
\over 
\vc{\pi}_{n+1}\vc{Q}_{n+1,n} 
\sum_{\ell=0}^n \vc{U}_{n,\ell}^*\vc{e}
}
=
{
\sum_{\ell\in\bbK} \vc{\pi}_{\ell}\vc{e}
\over 
\sum_{\ell=0}^n \vc{\pi}_{\ell}  \vc{e}
},
\qquad n \in \bbZ_{\ge K},
\label{eqn-r_n(wt{alpha}_n)}
\end{align}
where the second equality is due to (\ref{eqn:210618-01}). Combining (\ref{eqn:210712-01}) and (\ref{eqn-r_n(wt{alpha}_n)}) leads to
\begin{align*}
\sum_{\ell\in\bbK} \presub{(n)}\wh{\vc{\pi}}_{\ell}^* \vc{e}
\ge {
\sum_{\ell\in\bbK} \vc{\pi}_{\ell}\vc{e}
\over 
\sum_{\ell=0}^n \vc{\pi}_{\ell}  \vc{e}
}
\to \sum_{\ell\in\bbK} \vc{\pi}_{\ell}\vc{e}\quad \mbox{as $n \to \infty$}.
\end{align*}
Consequently, (\ref{lim-r_n(psi_n^*)}) holds. The proof has been completed. 

\medskip

Theorem~\ref{thm-E^*-01} shows that the new MIP-form solution (\ref{new-MIP-01}) is obtained by finding the optimal solution $\vc{\alpha}_n^*$ of Problem~\ref{prob-02}, that is, by finding an element in $\bbJ_n^*\subseteq \bbI_n^+$. This task can be lightened by the following theorem.
\begin{THM}\label{lem-J_n}
The set $\bbJ_n^*$ is given by
\begin{align}
\bbJ_n^* 
= \argmax_{j\in\bbI_n^+ \cap \bbO_n^+} {u_{n,\bbK}^*(j) \over u_n^*(j)},
\qquad n \in \bbZ_{\ge K},
\label{defn-J_n^*-02}
\end{align}
where $\bbO_n^+$ is defined as
\begin{equation}
\bbO_n^+= 
\{ i \in \bbM_n: (\vc{U}_n^*\vc{Q}_{n,n-1}\vc{e})_i > 0 \},
\qquad n \in \bbN.
\label{defn-J_n}
\end{equation}
\end{THM}

\medskip

\begin{proof}
Lemma~\ref{lem-optima-solution-LFP} states that, for any $ j \in \bbJ_n^*$, $u_{n,\bbK}^*(j)/u_{n}^*(j) = r_n(\vc{\alpha}_n^*) > 0$. This implies that $\bbJ_n^*$ does not include $j$ such that $u_{n,\bbK}^*(j) = 0$. It thus suffices to show that
\begin{equation}
j  \in \bbM_n\setminus\bbO_n^+ ~\Longrightarrow~ 
u_{n,\bbK}^*(j) = 0.
\label{relation-J_n^*}
\end{equation}
It follows from (\ref{defn-u_{n,N}^*}) and (\ref{eqn-U_{n,k}^*}) that
\begin{align}
\vc{u}_{n,\bbK}^*
&= \vc{U}_n^*\vc{Q}_{n,n-1} 
\cdot \sum_{\ell\in\bbK \cap [0,n-1]}
\!\! \!\! \!\! \vc{U}_{n-1,\ell}^* \vc{e},
\qquad n \in \bbZ_{\ge K}.
\label{recursion:u_{n,K}^*-org}
\end{align}
It also follows from (\ref{defn-J_n}) that
\begin{alignat*}{2}
\vc{U}_n^*\vc{Q}_{n,n-1}
&= 
\bordermatrix
{
         & \bbM_{n-1}
\cr
\bbO_n^+ & \mbox{\Large $\bm{*}$}
\cr
\bbM_n\setminus\bbO_n^+ & \vc{O}
},\qquad n \in\bbN.
\end{alignat*}
These two equations show that (\ref{relation-J_n^*}) holds, which completes the proof.
\end{proof}

\section{Advantages of the new MIP-form solution}\label{sec-discussion}

This section explains the advantages of our new MIP-form solution (\ref{new-MIP-01}) over the existing one (\ref{eqn-MIP-Masu19-QUESTA}). We first point out the drawbacks of the existing MIP-form solution. We then show that our new solution (\ref{new-MIP-01}) is quasi-algorithmically constructible (see Notion~\ref{defn-construction}) unlike the existing one. 

The existing MIP-form solution (\ref{eqn-MIP-Masu19-QUESTA}) has three drawbacks, compared with our new one. 
First, the existing solution (\ref{eqn-MIP-Masu19-QUESTA}) requires a parameter set $(\vc{v},b,\bbC) \in \cal{V}$ of Condition~\ref{cond-01} such that $\vc{v}$ satisfies Condition~\ref{cond-02} (see Proposition~\ref{prop-(n)wh{pi}^{dag}}). In some cases, Condition~\ref{cond-02} does not hold, and such an {\it unfavorable} case is provided in Appendix~\ref{append-cond-02}. Second, finding a parameter set $(b,\bbC,\vc{v}) \in \cal{V}$ may not be easy. Finally, it is of infinite computational complexity to obtain the optimal solution $\vc{\alpha}_n^{\dag}$ for Problem~\ref{prob-01}. This is because the objective function of Problem~\ref{prob-01} includes the infinite sums $\sum_{\ell=n+1}^{\infty} \vc{Q}_{k,\ell}\vc{v}_{\ell}$, $k \in \bbZ_{[0,n]}$ in (\ref{defn-y_n}).

Owing to these drawbacks, especially, to the last one, the existing MIP-form solution (\ref{eqn-MIP-Masu19-QUESTA}) for $\vc{\pi}$ is not quasi-algorithmically constructible. Nevertheless, constructing $\vc{\pi}$ via the existing solution (\ref{eqn-MIP-Masu19-QUESTA}) is {\it formally} expressed in an algorithm style equipped with a stopping criterion, as the {\it sequential update algorithm} \cite[Algorithm~1]{Masu19-QUESTA}. This is described in Algorithm~\ref{algo-01} below.
\begin{algorithm}[H]
\renewcommand{\labelenumi}{\arabic{enumi}.}
\renewcommand{\labelenumii}{(\alph{enumii})}
\renewcommand{\labelenumiii}{\roman{enumiii}.}
\caption{Sequential update algorithm for the existing MIP-form solution (\ref{eqn-MIP-Masu19-QUESTA})}\label{algo-01}
{\bf Input}: 
\begin{minipage}[t]{0.8\textwidth}
$\vc{Q}$, $\varepsilon \in (0,1)$, increasing sequence $\{n_{\ell} \in \bbN;\ell\in \bbZ_+\}$, and $(\vc{v},b,\bbC) \in \cal{V}$ 

\vspace{-2mm}

satisfying Conditions~\ref{cond-01} and \ref{cond-02}.
\end{minipage}

\smallskip

{\bf Output}: 
$\presub{(n)}\wh{\vc{\pi}}^{\dag} = (\presub{(n)}\wh{\vc{\pi}}_{k}^{\dag})_{k \in \bbZ_{[0,n]}}$,
where $n \in \bbZ_+$ is fixed when the iteration stops.

\vspace{-3mm}

\begin{enumerate}
\setlength{\parskip}{-0.0mm} 
\setlength{\itemsep}{-0.0mm} 
%
%
\item Set $n = 0$ and $\ell=1$.
\item Compute $\vc{U}_0^*$ by (\ref{eqn-U_{0,0}^*}), and compute $\vc{u}_0^*$ by (\ref{eqn-u_0^*}).
\item Iterate (a)--(d) below:

\vspace{-2mm}

\begin{enumerate}
\setlength{\parskip}{-0.0mm} 
\setlength{\itemsep}{-0.0mm} 
\item Increment $n$ by one.
\item Compute $\vc{U}_n^*=\vc{U}_{n,n}^*$ by (\ref{eqn-U_n^*}).
\item Compute $\vc{U}_{n,k}^*$, $k=0,1,\dots,n-1$, by (\ref{eqn-U_{n,k}^*}); and  compute $\vc{u}^*_n$ by (\ref{eqn-u_n^*}).
\item  If $n=n_{\ell}$, then perform the following:
\begin{enumerate}
\item Compute $\vc{y}_n$ by (\ref{defn-y_n}), and find $j_n^{\dag}$ satisfying (\ref{defn-j_n^{dag}}).
\item Compute $\presub{(n)}\wh{\vc{\pi}}_{k}^{\dag}$, $k=0,1,\dots,n$, by (\ref{defn-(n)wh{pi}_k^{dag}}).
\item If $\| \presub{(n_{\ell})}\wh{\vc{\pi}}^{\dag} - 
\presub{(n_{\ell-1})}\wh{\vc{\pi}}^{\dag}\|_1 < \varepsilon$, then stop the iteration; otherwise increment $\ell$ by one and return to step~(a).
\end{enumerate}
\end{enumerate}
\end{enumerate}
\end{algorithm}

\begin{REM}\label{rem-input}
In the original description of Algorithm~\ref{algo-01}, the following operation is inserted as the first step:
``Find $\vc{v}>\vc{0}$, $b>0$, and $\bbC \in \bbS$ such that Conditions~\ref{cond-01} and \ref{cond-02} hold" (see \cite[Algorithm~1]{Masu19-QUESTA}). However, in general, this cannot be performed algorithmically. Rather, $\vc{v}>\vc{0}$, $b>0$, and $\bbC \in \bbS$ are {\it input parameters} of the algorithm, as described in Algorithm~\ref{algo-01} above.
\end{REM}

\begin{REM}
In general, we have to truncate the infinite sums $\sum_{\ell=n+1}^{\infty} \vc{Q}_{k,\ell}\vc{v}_{\ell}$, $k \in \bbZ_{[0,n]}$ in order to compute $\vc{y}_n$ (see Remark~\ref{rem-LD-QBD} for an exceptional case).
\end{REM}

In contrast, our new MIP-form solution (\ref{new-MIP-01}) is quasi-algorithmically constructible. The new solution uses the sequence of probability vectors $\{\vc{\alpha}_n^*\}$, and each of them is an optimal solution for Problem~\ref{prob-02}. This problem is formulated without finding such a parameter set $(\vc{v},b,\bbC) \in \cal{V}$ satisfying Condition~\ref{cond-02}. Moreover, the objective function of Problem~\ref{prob-02} consists of only $\vc{u}_n^*$ and $\vc{u}_{n,\bbK}^*$ of finite computational complexity. Indeed, the finite computational complexity of $\vc{u}_n^*$ is guaranteed by the recursion (\ref{recursion-u_n^*}), and that of $\vc{u}_{n,\bbK}^*$ is guaranteed by the following recursion:
\begin{equation}
\vc{u}_{n,\bbK}^* 
= 
\left\{
\begin{array}{ll}
\dm\sum_{\ell\in\bbK} \vc{U}_{K,\ell}^*\vc{e}, & \quad n=K,
\\
\rule{0mm}{4mm}
\vc{U}_n^*\vc{Q}_{n,n-1} \cdot \vc{u}_{n-1,\bbK}^*,
&  \quad n \in \bbZ_{\ge K+1},
\end{array}
\right.
\label{recursion-u_{n,N}^*}
\end{equation}
which follows from (\ref{defn-u_{n,N}^*}) and (\ref{recursion:u_{n,K}^*-org}).  Therefore, it is of {\it finite} computational complexity to obtain the optimal solution $\vc{\alpha}_n^*$ for Problem~\ref{prob-02}. In summary, constructing our new MIP-form solution (\ref{new-MIP-01}) requires $\{\vc{\alpha}_n^*;n\in\bbZ_{\ge K}\}$, $\{\vc{u}_n^*;n\in\bbZ_+\}$, and $\{\vc{U}_{n,k}^*;n\in\bbZ_+,k\in\bbZ_{[0,n]}\}$, which are computed by iterating the recursive procedure consisting of an increasing but finite number of operations per iteration. Consequently, our new solution (\ref{new-MIP-01}) for $\vc{\pi}$ is quasi-algorithmically constructible.

Shown below is an algorithm-style description of computing our new MIP-form solution (\ref{new-MIP-01}) with a stopping criterion. This is the {\it sequential update algorithm} for our new solution.
\begin{algorithm}[H]
\renewcommand{\labelenumi}{\arabic{enumi}.}
\renewcommand{\labelenumii}{(\alph{enumii})}
\renewcommand{\labelenumiii}{\roman{enumiii}.}
\caption{Sequential update algorithm for the new MIP-form solution (\ref{new-MIP-01})}\label{algo-02}
{\bf Input}: $\vc{Q}$, $\varepsilon \in (0,1)$, finite set $\bbK \subset \bbZ_+$, and increasing sequence $\{n_{\ell} \in \bbZ_{\ge K};\ell\in \bbZ_+\}$.
\\ 
{\bf Output}: 
$\presub{(n)}\wh{\vc{\pi}}^* = (\presub{(n)}\wh{\vc{\pi}}_{k}^*)_{k \in {[0,n]}}$,
where $n \in \bbZ_{\ge K}$ is fixed when the iteration stops.

\vspace{-2mm}

\begin{enumerate}
\setlength{\parskip}{-0.0mm} 
\setlength{\itemsep}{-0.0mm} 
\item Set $n = 0$ and $\ell=1$.
\item Compute $\vc{U}_0^*$ by (\ref{eqn-U_{0,0}^*}), and compute $\vc{u}_0^*$ by (\ref{eqn-u_0^*}).
\item Iterate (a)--(d) below:

\vspace{-2mm}

\begin{enumerate}
\setlength{\parskip}{-0.0mm} 
\setlength{\itemsep}{-0.0mm} 
\item Increment $n$ by one.
\item Compute $\vc{U}_n^*=\vc{U}_{n,n}^*$ by (\ref{eqn-U_n^*}).
\item Compute $\vc{U}_{n,k}^*$, $k=0,1,\dots,n-1$, by (\ref{eqn-U_{n,k}^*}); and compute $\vc{u}^*_n$ by (\ref{eqn-u_n^*}).
\item  If $n=n_{\ell}\in\bbZ_{\ge K}$, then perform the following:
\begin{enumerate}
\item Compute $\vc{u}_{n,\bbK}^*$ by (\ref{recursion-u_{n,N}^*}), and find an element $j_n^*$ of $\bbJ_n^*$ given in (\ref{defn-J_n^*-02}).
\item Compute $\presub{(n)}\wh{\vc{\pi}}_{k}^*$, $k=0,1,\dots,n$, by (\ref{defn-(n)wh{pi}_k^*}).
\item If $\| \presub{(n_{\ell})}\wh{\vc{\pi}}^* - 
\presub{(n_{\ell-1})}\wh{\vc{\pi}}^*\|_1 < \varepsilon$, then stop the iteration; otherwise increment $\ell$ by one and return to step~(a).
\end{enumerate}
\end{enumerate}
\end{enumerate}
\end{algorithm}

\begin{REM}\label{rem-O_n^+}
To compute the probability vector $\presub{(n)}\wh{\vc{\pi}}^*=(\presub{(n)}\wh{\vc{\pi}}_0^*,\presub{(n)}\wh{\vc{\pi}}_1^*,\dots,\presub{(n)}\wh{\vc{\pi}}_n^*)$, Step (3.d.i) performs finding an element $j_n^*$ of $\bbJ_n^*$.
Finding an element $j_n^*$ of $\bbJ_n^*$ is equivalent to solving Problem~\ref{prob-02}, and it requires constructing the sets $\bbI_n^+$ and $\bbO_n^+$. This construction can be easily performed based on the probabilistic interpretation: The set $\bbI_n^+$ consists of the phases of $\bbL_n$ which accept direct (incoming) transitions from $\bbL_{n+1}$; the set $\bbO_n^+$ consists of the phases of $\bbL_n$ which are the starting points of (outgoing) paths leaving $\bbL_n$ eventually and reaching $\bbL_{n-1}$ while avoiding $\ol{\bbS}_n=\cup_{k=n+1}^{\infty}\bbL_k$. Note that the set $\bbO_n^+$ is obtained by identifying the positive rows of $\vc{U}_n^*\vc{Q}_{n,n-1}$, which is a crucial component of the recursion (\ref{eqn-U_{n,k}^*}) of $\{\vc{U}_{n,k}^*\}$. 
\end{REM}

\begin{REM}\label{rem-stopping-01}
The choice of finite set $\bbK\subset \bbZ_+$ can impact on the convergence speed of Algorithm~\ref{algo-02}. However, it would be difficult to discuss theoretically an optimal set $\bbK$. A reasonable choice of $\bbK$ would be $\bbZ_{[0,K]} = \{0,1,\dots,K\}$.
\end{REM}

\begin{REM}
A single run of Algorithm~\ref{algo-02} generates the probability vector $\presub{(n_{\ell})}\wh{\vc{\pi}}^*=(\presub{(n_{\ell})}\wh{\vc{\pi}}_0^*,\presub{(n_{\ell})}\wh{\vc{\pi}}_1^*,\dots,\presub{(n_{\ell})}\wh{\vc{\pi}}_{n_{\ell}}^*)$ with $\ell$ being equal to the value fixed on completion of all the operations. Recall that $\presub{(n_{\ell})}\wh{\vc{\pi}}^*$ is an LBC-LA truncation approximation to the stationary distribution vector $\vc{\pi}$. Thus, $\sum_{k=1}^{n_{\ell}}k^m\presub{(n_{\ell})}\wh{\vc{\pi}}_k^*$, $m\in\bbN$ can be considered an approximation to the moment vector $\sum_{k=1}^{\infty}k^m\vc{\pi}_k$ of the stationary distribution, though, as the order $m$ is larger, a smaller $\varep > 0$ would need to be chosen (as the parameter of the stopping criterion) in order to achieve a satisfactory accuracy of the approximation.
\end{REM}

Finally, in the two subsequent paragraphs, we mention the advantages of Algorithm~\ref{algo-02} based on our new MIP-form solution (\ref{new-MIP-01}), comparing with the existing algorithms \cite{M.Kimu18,Klim06,LiQuan05-RG,Shin98,Taki16} for UBH-MCs. These existing algorithms are classified into two types: (i) approximation by conditional stationary distribution \cite{M.Kimu18,Taki16}; (ii) approximation by level-homogenizing \cite{Klim06,LiQuan05-RG,Shin98}. Both types are originally designed to compute a single approximation to the stationary distribution, but performing them iteratively can generate a convergent sequence of approximations to the stationary distribution. Considering such applications of the existing algorithms, we compare them with our algorithm in the following.

The algorithms \cite{M.Kimu18,Taki16} compute a {\it single} conditional stationary distribution $\vc{\pi}^{(n)}$, as an approximation to $\vc{\pi}$, in a {\it single} run. To generate each approximation $\vc{\pi}^{(n)}$, the algorithms \cite{M.Kimu18,Taki16} have to perform a large (theoretically infinite) amount of {\it preliminary} computations involving the levels above $n$. The preliminary numerical results can be used to construct {\it better} approximations $\vc{\pi}^{(N)}$ with $N > n$, but to do this actually with satisfactory accuracy requires {\it further} preliminary computations involving {\it higher} levels. In contrast, Algorithm~\ref{algo-02} computes each tentative solution $\presub{(n)}\wh{\vc{\pi}}^*$ without any information on the higher levels above $n$ except $\vc{Q}_{n+1,n}$.

The other algorithms \cite{Klim06,LiQuan05-RG,Shin98} conduct level-homogenizing the transition law above some level, say level $n$, in order to compute an approximation to $\vc{\pi}$, denoted by $\vc{\pi}(n)$. More specifically, the approximation $\vc{\pi}(n)$ requires, as the input for computing itself, a (level-independent) $G$-matrix obtained by level-homogenizing the transition law above level $n$. Thus, each approximation $\vc{\pi}(n)$ is computed {\it from scratch}, which is time-consuming. In contrast, Algorithm~\ref{algo-02} {\it sequentially} constructs tentative solutions $\{\presub{(n)}\wh{\vc{\pi}}^*\}$ convergent to the stationary distribution vector $\vc{\pi}$, making the most of the components of the ``old" tentative solutions.

\section{Special cases free from solving Problem~\ref{prob-02}}\label{sec-concluding}

This section discusses special cases where an MIP-form solution for $\vc{\pi}$ is constructed without solving Problem~\ref{prob-02}. We first show a basic theorem on this matter, and then derive its corollary that implies what cases are free from solving Problem~\ref{prob-02}. Finally, we mention representative examples of such favorable cases.

The following theorem provides the basis for our discussion here.
\begin{THM}\label{Prop:MIP} 
Let $N$ be a nonnegative integer, and let $\vc{\varphi}_n$, $n \in \bbZ_{\ge N+1}$ denote a $1 \times M_n$ probability vector such that
\begin{equation}
\left| 
\vc{\varphi}_n
-
{ \vc{\pi}_n \over \vc{\pi}_n\vc{e} }
\right| \le \varep_n { \vc{\pi}_n \over \vc{\pi}_n\vc{e} },
\qquad n \in \bbZ_{\ge N+1},
\label{cond:varphi_n}
\end{equation}
where $\{\varep_n \in [0,1);n \in \bbZ_{\ge N+1}\}$ is a sequence such that $\lim_{n\to\infty}\varep_n = 0$. Furthermore, for $n \in \bbZ_{\ge N}$, let $\presub{(n)}\bv{\vc{\pi}}:=(\presub{(n)}\bv{\vc{\pi}}_0,\presub{(n)}\bv{\vc{\pi}}_1,\dots,\presub{(n)}\bv{\vc{\pi}}_n)$ denote a probability vector such that
\begin{align}
\presub{(n)}\bv{\vc{\pi}}_k
= 
{ 
\vc{\varphi}_{n+1} \prod_{m=k}^{n} \downarrow \vc{U}_m
\over 
\vc{\varphi}_{n+1} 
\sum_{\ell=0}^{n} (\prod_{m=\ell}^{n} \downarrow \vc{U}_m) \vc{e}
},\qquad k\in\bbZ_{[0,n]}.
\label{defn:(n)bv{pi}_k}
\end{align}
We then have
\begin{align}
\left\|
\presub{(n)}\bv{\vc{\pi}} - \vc{\pi}
\right\|_1
\le {2 \varep_{n+1} \over 1 - \varep_{n+1}}
+ 2\sum_{k=n+1}^{\infty} \vc{\pi}_k\vc{e}, \qquad n \in \bbZ_{\ge N}.
\label{new-MIP-form-02}
\end{align}
\end{THM}

\begin{proof}
Fix $n \in \bbZ_{\ge N}$ arbitrarily, and note that (\ref{new-MIP-form-02}) holds if
\begin{align}
\left|
\presub{(n)}\bv{\vc{\pi}}_k
-
{
\vc{\pi}_k
\over 
\sum_{\ell=0}^n \vc{\pi}_{\ell}\vc{e}
}
\right|
\le
{
2\varep_{n+1}
\over
1 - \varep_{n+1}
}
{
\vc{\pi}_k
\over 
\sum_{\ell=0}^n \vc{\pi}_{\ell}\vc{e}
},\qquad k \in \bbZ_{[0,n]}.
\label{eqn:210622-02}
\end{align}
Indeed, the triangle inequality and (\ref{defn-pi^{(N)}}) yield 
\begin{align*}
\left\| \presub{(n)}\bv{\vc{\pi}} - \vc{\pi} \right\|_1
&\le 
\left\| \presub{(n)}\bv{\vc{\pi}} - \vc{\pi}^{(n)} \right\|_1 
+
\left\| \vc{\pi}^{(n)} - \vc{\pi} \right\|_1
\nonumber
\\
&= \left\| \presub{(n)}\bv{\vc{\pi}} - \vc{\pi}^{(n)} \right\|_1
+ 2\sum_{k=n+1}^{\infty} \vc{\pi}_k\vc{e}
\nonumber
\\
&= \sum_{k=0}^n 
\left| 
\presub{(n)}\bv{\vc{\pi}}_k 
- { \vc{\pi}_k \over \sum_{\ell=0}^{n} \vc{\pi}_{\ell}\vc{e}} 
\right| \vc{e}
+ 2\sum_{k=n+1}^{\infty} \vc{\pi}_k\vc{e},
\end{align*}
and then, combining this and (\ref{eqn:210622-02}) results in (\ref{new-MIP-form-02}).

To complete the proof, we show that (\ref{eqn:210622-02}) holds. Equation~(\ref{cond:varphi_n}) implies that
\begin{align}
(1 - \varep_{n+1})
{\vc{\pi}_{n+1} \over \vc{\pi}_{n+1}\vc{e}}
\le
\vc{\varphi}_{n+1} 
\le
(1 + \varep_{n+1})
{ \vc{\pi}_{n+1} \over \vc{\pi}_{n+1}\vc{e}}.
\label{eqn:210622-01}
\end{align}
Applying (\ref{eqn:210622-01}) to (\ref{defn:(n)bv{pi}_k}) yields
\begin{align*}
{
1 - \varep_{n+1}
\over
1 + \varep_{n+1}
}
{
\vc{\pi}_{n+1} \prod_{m=k}^{n} \downarrow \vc{U}_m
\over 
\vc{\pi}_{n+1} \sum_{\ell=0}^n (\prod_{m=\ell}^{n} \downarrow \vc{U}_m)\vc{e}
}
\le
\presub{(n)}\bv{\vc{\pi}}_k
\le
{
1 + \varep_{n+1}
\over
1 - \varep_{n+1}
}
{
\vc{\pi}_{n+1} \prod_{m=k}^{n} \downarrow \vc{U}_m
\over 
\vc{\pi}_{n+1} \sum_{\ell=0}^n (\prod_{m=\ell}^{n} \downarrow \vc{U}_m)\vc{e}
}.
\end{align*}
Furthermore, substituting (\ref{eqn:210618-01}) into the above inequality leads to
\begin{align*}
{
1 - \varep_{n+1}
\over
1 + \varep_{n+1}
}
{
\vc{\pi}_k
\over 
\sum_{\ell=0}^n \vc{\pi}_{\ell}\vc{e}
}
\le
\presub{(n)}\bv{\vc{\pi}}_k
\le
{
1 + \varep_{n+1}
\over
1 - \varep_{n+1}
}
{
\vc{\pi}_k
\over 
\sum_{\ell=0}^n \vc{\pi}_{\ell}\vc{e}
},\qquad k \in \bbZ_{[0,n]},
\end{align*}
and thus
\begin{align*}
\left|
\presub{(n)}\bv{\vc{\pi}}_k
- 
{
\vc{\pi}_k
\over 
\sum_{\ell=0}^n \vc{\pi}_{\ell}\vc{e}
}
\right|
\le 
{
2\varep_{n+1}
\over
1 - \varep_{n+1}
}
{
\vc{\pi}_k
\over 
\sum_{\ell=0}^n \vc{\pi}_{\ell}\vc{e}
},\qquad k \in \bbZ_{[0,n]},
\end{align*}
which shows that (\ref{eqn:210622-02}) holds. The proof has been completed.
\end{proof}

\begin{REM}
Substituting 
\begin{align*}
\vc{\alpha}_n 
= 
{
\vc{\varphi}_{n+1}\vc{Q}_{n+1,n} \over \vc{\varphi}_{n+1}\vc{Q}_{n+1,n}\vc{e}
}
\end{align*}
into the matrix product form (\ref{eqn-(n_s)pi_{n_s,l}}) of the LBC-LA truncation approximation $\presub{(n)}\wh{\vc{\pi}}=(\presub{(n)}\wh{\vc{\pi}}_0,\presub{(n)}\wh{\vc{\pi}}_1,\dots,\presub{(n)}\wh{\vc{\pi}}_n)$, and using (\ref{defn:U_m}), we obtain $\presub{(n)}\wh{\vc{\pi}}=\presub{(n)}\bv{\vc{\pi}}$. Therefore, $\presub{(n)}\bv{\vc{\pi}}$ is an LBC-LA truncation approximation to $\vc{\pi}$.
\end{REM}

Theorem~\ref{Prop:MIP} implies in what cases we do not have to solve Problem~\ref{prob-02}. However, Theorem~\ref{Prop:MIP} requires us to identify $\vc{\pi}_n /(\vc{\pi}_n\vc{e})$ to some extent, which is not easy in general. We thus provide a more effective (but more restrictive) result that can actually contribute to the construction of MIP-form solutions without solving Problem~\ref{prob-02}.  
\begin{COR}\label{cor:special-case-02}
Suppose that there exists some $N \in \bbZ_+$ such that $\bbM_n = \bbM:=\{1,2,\dots,M\} \subset \bbN$ for all $n \in \bbZ_{\ge N+1}$. Furthermore, there exists some $1 \times M$ probability vector $\vc{\varpi}$ such that
\begin{equation}
\lim_{n \to \infty} 
{\vc{\pi}_n \over \vc{\pi}_n\vc{e}} = \vc{\varpi} > \vc{0}.
\label{lim-pi_n/pi_n*e}
\end{equation}
We then have
\begin{eqnarray*}
\lim_{n\to\infty}
\left\|
\vc{\pi}
-
\left(
{ 
\vc{\varpi} \prod_{m=k}^{n} \downarrow \vc{U}_m
\over 
\vc{\varpi} 
\sum_{\ell=0}^{n} (\prod_{m=\ell}^{n} \downarrow \vc{U}_m) \vc{e}
}
\right)_{\!\!k\in\bbZ_{[0,n]}}
\right\|_1 = 0.
\end{eqnarray*}
\end{COR}

\begin{proof}
According to (\ref{lim-pi_n/pi_n*e}) and Theorem~\ref{Prop:MIP}, it suffices to show that (\ref{cond:varphi_n}) holds with $\vc{\varphi}_n = \vc{\varpi}$ for all $n \in \bbZ_{\ge N+1}$. Equation~(\ref{lim-pi_n/pi_n*e}) implies that there exists some $\{\delta_n \in [0,1);n \in \bbZ_{\ge N+1}\}$ such that $\lim_{n\to\infty}\delta_n=0$ and
\begin{align*}
\left| {\vc{\pi}_n \over \vc{\pi}_n\vc{e}} - \vc{\varpi} \right|
\le \delta_n \vc{\varpi},\qquad n \in \bbZ_{\ge N+1}.
\end{align*}
This inequality is equivalent to
\begin{align}
{1 \over 1 + \delta_n}{\vc{\pi}_n \over \vc{\pi}_n\vc{e}}
\le
\vc{\varpi}
\le
{1 \over 1 - \delta_n}{\vc{\pi}_n \over \vc{\pi}_n\vc{e}},
\qquad n \in \bbZ_{\ge N+1}.
\label{eqn:210622-03}
\end{align}
We now fix
\begin{align*}
\varep_n = {1 \over 1 - \delta_n} - 1 = {\delta_n \over 1 - \delta_n},
\qquad n \in \bbZ_{\ge N+1}.
\end{align*}
We then have $\lim_{n\to\infty}\varep_n=0$ and 
\begin{alignat*}{2}
{1 \over 1 - \delta_n} 
&= 1 + \varep_n, & \qquad n &\in \bbZ_{\ge N+1},
\\
{1 \over 1 + \delta_n} 
&= {1 + \varep_n \over 1 + 2\varep_n}
\ge 1 - \varep_n, & \qquad n &\in \bbZ_{\ge N+1}.
\end{alignat*}
Thus, (\ref{eqn:210622-03}) yields
\begin{align*}
(1 - \varep_n){\vc{\pi}_n \over \vc{\pi}_n\vc{e}}
\le
\vc{\varpi}
\le
(1 + \varep_n){\vc{\pi}_n \over \vc{\pi}_n\vc{e}},
\qquad n \in \bbZ_{\ge N+1},
\end{align*}
which shows that (\ref{cond:varphi_n}) holds with $\vc{\varphi}_n = \vc{\varpi}$ for all $n \in \bbZ_{\ge N+1}$. 
\end{proof}

A typical example covered by Corollary~\ref{cor:special-case-02} is a UBH-MC with {\it block monotonicity} (see \cite[Definition~3.2]{Masu17-LAA}). If the generator $\vc{Q}$ is block monotone, then, for all $k \in \bbZ_+$, $\sum_{\ell=\max(k-1,0)}^{\infty}\vc{Q}_{k,\ell}$ is constant and thus is interpreted as the generator of the underlying Markov chain on the phase set $\bbM$, which controls the level process (see \cite[Lemma~3.1]{Masu17-LAA}). In this case, $\vc{\varpi}$ is equal to the stationary distribution vector of the underlying Markov chain. 

Another example is the (level-independent) M/G/1-type Markov chain. For this chain, several sufficient conditions are established for the tail asymptotics of $\vc{\pi}=(\vc{\pi}_0,\vc{\pi}_1,\dots)$ (see \cite{Kimu10,Kimu13,Kimu17,Masu16-ANOR}) such that, for some $c > 0$, $\vc{\mu} > \vc{0}$, and random variable $Y$ on $\bbZ_+$,
\begin{equation*}
\lim_{n\to\infty} {\sum_{k=n+1}^{\infty} \vc{\pi}_k \over \PP(Y>n)} 
= c\vc{\mu}.
\end{equation*}
Moreover, (assuming additional conditions, if required), we can obtain (see, e.g., \cite[Sections~4 and 5]{Kimu13})
\begin{equation*}
\lim_{n\to\infty} { \vc{\pi}_n \over \PP(Y=n)} = c\vc{\mu},
\end{equation*}
which implies that (\ref{lim-pi_n/pi_n*e}) holds.

\section{Concluding Remarks}\label{sec:concluding-remarks}

We have established a quasi-algorithmic construction of the exact (not approximate) and whole (not partial) stationary distribution vector $\vc{\pi}$ in upper block-Hessenberg Markov chains (UBH-MCs). The core of this theoretical construction is that all we have to do is to solve {\it just one} linear fractional programing (LFP) problem (Problem~\ref{prob-02}) in {\it each} iteration for the construction. We have also presented some special cases free from solving the LFP problem. To find other such special cases is an interesting future task. A possible case would be the class of UBH-MCs with {\it asymptotically block-Toeplitz} structure. The generator $\vc{Q}$ of such a chain satisfies (\ref{defn-Q}) and the following:
$M_n = M \in \bbN$ for all sufficiently large $n$ and
\begin{align*}
\lim_{n\to\infty}\vc{Q}_{n,n+k} = \vc{Q}_k,\qquad k \in \bbZ_{\ge-1},
\end{align*}
where $\sum_{k=-1}^{\infty}\vc{Q}_k$ is an essential $Q$-matrix. In \cite{Klim06}, this special UBH-MC is referred to as a {\it multi-dimensional asymptotically quasi-Toeplitz Markov chain}.

\appendix

\section{An example violating Condition~\ref{cond-02}}\label{append-cond-02}

This section presents an example of UBH $Q$-matrices for which the existing MIP-form solution in \cite{Masu19-QUESTA} does not hold. To describe that example, we assume that $\vc{Q}$ is the generator of a specific M/G/1-type Markov chain with power-like level increments. For such a generator $\vc{Q}$, we show that Condition~\ref{cond-02} does not hold for any parameter set $(\vc{v},b,\bbC)$ of Condition~\ref{cond-01} and therefore the existing MIP-form solution is not established.

Let $\bbS = \bbZ_+ \times \bbM = \bbZ_+ \times \{1,2,\dots,M\}$, and suppose that $\vc{Q}=(q(k,i;\ell,j))_{(k,i;\ell,j)\in\bbS}$ is an M/G/1-type generator given by
\begin{align}
\vc{Q}
= \left(
\begin{array}{ccccc}
\vc{A}_{-1}+\vc{A}_0 &
\vc{A}_1 &
\vc{A}_2 &
\vc{A}_3 &
\cdots 
\\
\vc{A}_{-1} &
\vc{A}_0 &
\vc{A}_1 &
\vc{A}_2 &
\cdots 
\\
\vc{O} 	&
\vc{A}_{-1} &
\vc{A}_0 &
\vc{A}_1 &
\cdots 
\\
\vc{O} 	&
\vc{O} 	&
\vc{A}_{-1} &
\vc{A}_0 &
\cdots 
\\
\vdots	&
\vdots	&
\vdots	&
\ddots	&
\ddots	
\end{array}
\right),
\label{defn:Q-MG1}
\end{align}
where $\vc{A}_k:=(A_{k,i.j})_{i,j\in\bbM}$, $k\in\bbZ_{\ge -1}$, is an $M \times M$ matrix. We then assume the following.  
\begin{ASSU}\label{assumpt-appendix-B}
(i) The generator $\vc{Q}$ in (\ref{defn:Q-MG1}) is irreducible; (ii) $\sum_{k=-1}^{\infty}\vc{A}_k$ is an essential $Q$-matrix with stationary distribution vector~$\vc{\varpi}$; (iii) $\vc{\varpi}\sum_{k=-1}^{\infty} k\vc{A}_k\vc{e} < 0$; and (iv) $\lim_{k\to\infty}k^3\vc{A}_k = \vc{C}_A$ for some nonnegative matrix $\vc{C}_A \neq \vc{O}$.
\end{ASSU}

We confirm that there exists some $(\vc{v},b,\bbC) \in \cal{V}$ of Condition~\ref{cond-01} under Assumption~\ref{assumpt-appendix-B}. Uniformizing  the generator $\vc{Q}$ (see, e.g., \cite[Problem~II.4.1]{Asmu03}) yields an M/G/1 type stochastic matrix that satisfies the stability condition \cite[Chapter~XI, Proposition~3.1]{Asmu03} of GI/G/1-type Markov chains (including  M/G/1-type ones) due to the conditions (i)--(iii) of Assumption~\ref{assumpt-appendix-B}. Therefore, Condition~\ref{cond-01} holds for the generator $\vc{Q}$ satisfying Assumption~\ref{assumpt-appendix-B} (see, e.g., \cite[Theorem~1.1]{Kont16}). 

We confirm that Condition~\ref{cond-02} does not hold for any $(\vc{v},b,\bbC) \in \cal{V}$ of Condition~\ref{cond-01} under Assumption~\ref{assumpt-appendix-B}. Fix $(\vc{v},b,\bbC) \in \cal{V}$ arbitrarily. Since $\bbC$ is a finite subset of the state space $\bbS = \bbZ_+ \times \bbM$, there exists some $N \in \bbZ_+$ such that $\bbC \subseteq \bbS_N=\bbZ_{[0,N]} \times \bbM$. Thus, fix $N \in \bbZ_+$ as such, and let $\{(X_t,J_t); t \ge 0\}$ denote an M/G/1-type Markov chain on state space $\bbS$ with the generator $\vc{Q}$ satisfying Assumption~\ref{assumpt-appendix-B}. It then follows from \cite[Theorem~1.1]{Kont16} that there exists some $\gamma > 0$ such that 
\begin{align}
\gamma (v(k,i) + 1)
&\ge \EE[T_{\bbC}(1) \mid (X_0,J_0) = (k, i)]
\nonumber
\\
&\ge \EE[T_{\bbC} \mid (X_0,J_0) = (k, i)]
\quad \mbox{for all $(k,i) \in \bbS$},
\label{eqn:200703-01}
\end{align}
where $T_{\bbC}(1) = \inf\{t \ge 1: (X_t,J_t) \in \bbC\}$ and $T_{\bbC} = \inf\{t > 0: (X_t,J_t) \in \bbC\}$. It also follows from $\bbC \subseteq \bbS_N$ and the M/G/1-type structure of $\vc{Q}$ that there exists some $\delta > 0$ such that
\begin{equation}
\EE[T_{\bbC} \mid (X_0,J_0) = (k, i)] \ge  \delta k
\quad \mbox{for all $(k,i) \in \ol{\bbS}_N$}.
\label{eqn:200703-02}
\end{equation}
Combining (\ref{eqn:200703-01}) and (\ref{eqn:200703-02}), we have
\begin{align}
v(k,i) \ge {\delta \over \gamma} k - 1 \ge 0 \quad \mbox{for all $k \ge N':=\max\left(N+1,\dm{\gamma \over \delta} \right)$ and $i \in \bbM$}.
\label{eqn:210721-01}
\end{align}
Furthermore, from \cite[Theorem~4.1.1]{Kimu13} and Assumption~\ref{assumpt-appendix-B} (especially, the condition~(iv)), we have
\begin{equation}
\lim_{k\to\infty}k^2 \vc{\pi}_k 
= c \vc{\varpi},
\label{cond-tai-pi}
\end{equation}
for some constant $c > 0$. Using (\ref{defn:Q-MG1}), (\ref{eqn:210721-01}), and (\ref{cond-tai-pi}), we obtain
\begin{align*}
\vc{\pi}|\vc{Q}|\vc{v} 
& \ge
\nonumber \sum_{(k,i) \in \bbS} \pi(k,i) |q(k,i;k,i)| v(k,i)
\\
&\ge \sum_{k=N'}^{\infty} \sum_{i \in\bbM} \pi(k,i) |A_{0,i,i}| v(k,i)
\nonumber
\\
&\ge \min_{j\in\bbM}|A_{0,j,j}| \cdot
\sum_{k=N'}^{\infty} 
\left( {\delta \over \gamma} k - 1 \right) \vc{\pi}_k\vc{e}=\infty.
\end{align*}
Therefore, Condition~\ref{cond-02} does not hold for any $(\vc{v},b,\bbC) \in \cal{V}$ in the present example.

\section{Convergent approximations to the stationary distribution vector of an   essential $Q$-matrix}

The purpose of this section is to provide the fundamental results needed in the proof of the main theorem (Theorem~\ref{thm-E^*-01}). This section consists of two subsections. Section~\ref{subsec:defn-peroper-Q} describes several notions associated with essential $Q$-matrices. Section~\ref{subsec:truncation} presents a necessary and sufficient condition for the convergence of a sequence of approximations to the stationary distribution vector of an essential $Q$-matrix.

\subsection{Definition of basic notions}\label{subsec:defn-peroper-Q}

This subsection introduces the notions used in the next subsection. We first redefine the symbols $\bbS$ and $\vc{Q}$ introduced in the body of this paper. We next define several notions, such as essential $Q$-matrices, subinvariant measures, irreducibility, and recurrence. We then present a proposition on the recurrence and subinvariant measures of an irreducible essential $Q$-matrix.

Let $\bbS$ denote an arbitrary countable set, and let $\vc{Q}:=(q(i,j))_{i,j\in\bbS}$ denote an arbitrary {\it $Q$-matrix} (see, e.g., \cite[page 64]{Ande91}), that is, a diagonally dominant matrix such that
\begin{alignat*}{2}
-\infty \le q(i,i) &\le 0, &\qquad  i &\in \bbS,
\\ 
0 \le q(i,j) &<\infty, &\qquad  i &\in \bbS,~j \in \bbS \setminus\{i\},
\\
\sum_{j \in \bbS\setminus\{i\}} q(i,j) &\le -q(i,i), & i &\in \bbS.
\end{alignat*}

The following is the definition of essential $Q$-matrices.
\begin{DEF}\label{defn-essential}
The $Q$-matrix $\vc{Q}$ is said to be {\it essential} if and only if the following hold for all $i \in \bbS$: (i) $q(i,i)$ is finite; (ii) $\sum_{j \in \bbS} q(i,j) = 0$; (iii) $q(i,i) < 0$. Note that any essential $Q$-matrix can be interpreted as the infinitesimal generator of a continuous-time Markov chain. Thus, an essential $Q$-matrix may be referred to as an {\it infinitesimal generator} or {\it generator}, depending on the context.
\end{DEF}

\begin{REM}\label{rem:stable_conservative}
The conditions (i) and (ii), respectively, imply that $\vc{Q}$ is stable and conservative (see, e.g., \cite[Definition~13.3.10]{Brem20}). 
\end{REM}

\begin{REM}
An essential $Q$-matrix is not necessarily {\it regular (or equivalently, non-explosive)}. Indeed, provided that $\vc{Q}$ is an essential $Q$-matrix, $\vc{Q}$ is regular if and only if, for any $\lambda > 0$, the system of equations
\[
\vc{Q}\vc{x} = \lambda \vc{x}\quad \mbox{with $\vc{x}:=(x(i))_{i\in\bbS} \ge \vc{0}$}
\]
has no bounded solution other than $\vc{x}=\vc{0}$ (see, e.g., \cite[Chapter~2, Theorem~2.7]{Ande91} and \cite[Theorem~13.3.11]{Brem20}).
\end{REM}

To proceed further, we assume that $\vc{Q}$ is an essential $Q$-matrix, and then define $\{\Phi(t);t\ge0\}$ as a Markov chain on $\bbS$ having this essential $\vc{Q}$ as its generator. For later use, we also define $\vc{P}$ as a stochastic matrix such that
\begin{equation}
\vc{P} = \vc{I} + \diag\{-\vc{Q}\}^{-1}\vc{Q},
\label{defn-P}
\end{equation}
where $\diag\{-\vc{Q}\}$ is a diagonal matrix whose diagonal elements are identical to those of $(-\vc{Q})$. Note that $\vc{P}$ is interpreted as the transition probability matrix of an embedded discrete-time Markov chain for $\{\Phi(t)\}$ with generator $\vc{Q}$ (see \cite[Section~13.3.2]{Brem20}). Hence, we call $\vc{P}$ the {\it embedded transition probability matrix} of $\{\Phi(t)\}$.

For the essential $Q$-matrix $\vc{Q}$, we define subinvariant measures and stationary distribution vectors.
\begin{DEF}\label{defn-subinvariant}
Let $\vc{\mu}:=(\mu(j))_{j\in\bbS}$ denote an arbitrary nonnegative and nonzero vector.
\begin{enumerate}
\item If $\vc{\mu}\vc{Q} \le \vc{0}$, $\vc{\mu}$ is said to be a {\it subinvariant measure} of $\vc{Q}$ or the Markov chain $\{\Phi(t)\}$. Furthermore, if a subinvariant measure $\vc{\mu}$ satisfies $\vc{\mu}\vc{Q} = \vc{0}$, $\vc{\mu}$ is said to be {\it invariant or stationary}. 
\item If an invariant measure $\vc{\mu}$ satisfies $\vc{\mu}\vc{e} = 1$, $\vc{\mu}$ is said to be a {\it stationary distribution vector} of $\vc{Q}$ or the Markov chain $\{\Phi(t)\}$. 
\end{enumerate}
\end{DEF}

Finally, we provide the definition of irreducibility, recurrence, and their related notions, and present a proposition that summarizes basic results on the recurrence and subinvariant measures of an irreducible essential $Q$-matrix.
\begin{DEF}
{}\hfill
\begin{enumerate}
\item An essential $Q$-matrix $\vc{Q}$ and its Markov chain $\{\Phi(t)\}$ are {\it irreducible} (resp.\ {\it transient}, {\it recurrent}) if and only if the embedded transition probability matrix $\vc{P}$ in (\ref{defn-P}) is {\it irreducible} (resp.\ {\it transient}, {\it recurrent}) (see \cite[Definitions~13.4.1 and 13.4.2]{Brem20}). 
\item An irreducible essential $Q$-matrix $\vc{Q}$ and its Markov chain $\{\Phi(t)\}$ are {\it ergodic} (i.e., {\it positive recurrent}) if and only if there exists a stationary distribution vector, or equivalently, a summable invariant measure unique up to constant multiples (see \cite[Definition~13.4.8 and Theorem~13.4.10]{Brem20}). 
\end{enumerate}

\end{DEF}

\begin{PROP}\label{prop-subinvariant-measure}
Suppose that $\vc{Q}$ is an irreducible essential $Q$-matrix. The following statements hold:
\begin{enumerate}
\item There always exists a subinvariant measure of $\vc{Q}$. 
\item Any subinvariant measure of $\vc{Q}$ is positive, that is, its elements are all positive.
\item If $\vc{Q}$ is recurrent, then it has an invariant measure, which is unique up to constant multiples.
\item Any subinvariant measure of recurrent $\vc{Q}$ is, in fact, an invariant measure.
\item The matrix $\vc{Q}$ has a subinvariant measure that is not invariant if and only if it is transient.
\end{enumerate}
\end{PROP}

\begin{proof}
Let $\vc{\mu}=(\mu_i)_{i\in\bbS}$ be a nonnegative and nonzero vector, and let
\begin{equation}
\vc{\eta}:=(\eta_i)_{i\in\bbS} = \vc{\mu} \diag\{-\vc{Q}\} \ge \vc{0}, \neq \vc{0}.
\label{defn-eta}
\end{equation}
It then follows from (\ref{defn-P}) and (\ref{defn-eta}) that
\begin{eqnarray*}
\vc{\eta}\vc{P}
&=& \vc{\eta} + \vc{\mu}\vc{Q}.
\end{eqnarray*}
Thus, $\vc{\eta}$ is a subinvariant (resp.\ an invariant) measure of the embedded transition probability matrix $\vc{P}$ if and only if $\vc{\mu}\vc{Q} \le \vc{0}$ (resp.~$\vc{\mu}\vc{Q} = \vc{0}$) or equivalently, $\vc{\mu}$ is a subinvariant (resp.\ an invariant) measure of the $Q$-matrix $\vc{Q}$ (see \cite[Definition~5.3]{Sene06}). Furthermore, since the $Q$-matrix $\vc{Q}$ is essential (that is, $-q(i,i) > 0$ for all $i \in \bbS$), 
\begin{align*}
\mu_i > 0 \Longleftrightarrow \eta_i > 0.
\end{align*}
Therefore, the present proposition is an immediate consequence of \cite[Theorem~5.4]{Sene06} (together with its corollary) and \cite[Lemmas~5.5 and 5.6]{Sene06}. That implication is as follows:
\begin{align*}
\mbox{\cite[Lemma 5.5]{Sene06}} &\Longrightarrow
\mbox{Statement~(i)};
\\
\mbox{\cite[Lemma 5.6]{Sene06}} &\Longrightarrow
\mbox{Statement~(ii)};
\\
\mbox{\cite[Theorem 5.4]{Sene06}} &\Longrightarrow
\mbox{Statement~(iii)};
\\
\mbox{\cite[The corollary of Theorem 5.4]{Sene06}} &\Longrightarrow
\mbox{Statement~(iv)};
\\
\mbox{\cite[Theorem 5.4]{Sene06}} &\Longrightarrow
\mbox{Statement~(v)}.
\end{align*}
The proof has been completed. 
\end{proof}

\subsection{Convergence condition of approximations}\label{subsec:truncation}

This subsection proves a single theorem (Theorem~\ref{thm-appen-03}), which provides a necessary and sufficient condition for the convergence of a sequence of approximations to the stationary distribution of an essential $Q$-matrix. That necessary and sufficient condition is crucial to the proof of Theorem~\ref{thm-E^*-01}, the main theorem of this paper.

We begin with the following assumption.
\begin{ASSU}\label{assumpt-appendix}
\hfill
\begin{enumerate}
\item The essential $Q$-matrix $\vc{Q}$ is ergodic with the unique stationary distribution vector $\vc{\pi}:=(\pi(j))_{j\in\bbS} > \vc{0}$ (the uniqueness and positivity of $\vc{\pi}$ is due to Proposition~\ref{prop-subinvariant-measure}). 
\item For all $n \in \bbZ_+$, $\vc{Q}_n:=(q_n(i,j))_{i,j\in\bbS}$ is a stable and conservative $Q$-matrix (see Remark~\ref{rem:stable_conservative}) such that $\lim_{n\to\infty}q_n(i,j) = q(i,j)$ for all $i,j \in \bbS$.
\item Each $\vc{Q}_n$ has at least one stationary distribution vector, denoted by $\vc{\pi}_n:=(\pi_n(j))_{j\in\bbS}$, which can be considered an approximation to $\vc{\pi}$.
\end{enumerate}
\end{ASSU}

\begin{REM}
The vector $\vc{\pi}_n$ in Assumption~\ref{assumpt-appendix} is redefined here and thus is completely different from ``$\vc{\pi}_n$" used as the subvector of $\vc{\pi}$ in the body of this paper.
\end{REM}

Under Assumption~\ref{assumpt-appendix}, we have the next lemma, which is used to prove the theorem of this subsection.
\begin{LEM}\label{lem-appen-02}
Suppose that Assumption~\ref{assumpt-appendix} holds. Fix $j_0 \in \bbS$ arbitrarily, and let $\bbH$ denote an infinite subset of $\bbZ_+$ such that
$\{\pi_n(j_0);n\in\bbH\}$ converges to some $\alpha \in [0,\infty)$, that is,
\begin{equation}
\alpha
=\lim_{\scriptstyle n\to\infty \atop \scriptstyle n \in \bbH}\pi_n(j_0).
\label{defn-alpha}
\end{equation}
Under these conditions, we have the statements (i) and (ii) below:
\begin{enumerate}
\item It holds that
\begin{align}
\lim_{\scriptstyle n\to\infty \atop \scriptstyle n \in \bbH}\pi_n(j) 
&= {\alpha \over \pi(j_0)} \pi(j)
\quad\mbox{for all $j \in \bbS$},
\label{eqn-lim-pi_n(j)}
\\
\alpha &\le \pi(j_0).
\label{ineqn-alpha}
\end{align}
\item Furthermore, if $\alpha > 0$, then
\begin{equation}
\lim_{\scriptstyle n\to\infty \atop \scriptstyle n \in \bbH}
\left\|
{\vc{\pi}_n^{\bbA} \over \vc{\pi}_n^{\bbA} \vc{e} }
-
{\vc{\pi}^{\bbA} \over \vc{\pi}^{\bbA} \vc{e} }
\right\|_1
=0\quad
\mbox{for all finite $\bbA \subset \bbS$},
\label{eqn-lim-pi_n^A}
\end{equation}
where, for any vector $\vc{x}:=(x(j))_{j\in\bbS}$, $\vc{x}^{\bbA}:=(x^{\bbA}(j))_{j\in\bbS}$ denotes a vector such that
\[
x^{\bbA}(j) 
=
\left\{
\begin{array}{ll}
x(j), & \quad j \in \bbA,
\\
0,    & \quad j \not\in \bbA.
\end{array}
\right.
\]
\end{enumerate}
\end{LEM}

\begin{proof}
We begin with the proof of the statement (i). To prove (\ref{eqn-lim-pi_n(j)}) and (\ref{ineqn-alpha}), it suffices to show that, for all infinite $\bbH' \subseteq \bbH$,
\begin{eqnarray}
\liminf_{\scriptstyle n\to\infty \atop \scriptstyle n \in \bbH'}
\pi_n(j)
&=& {\alpha \over \pi(j_0)} \pi(j)
\le  \pi(j),\qquad j \in \bbS.
\label{eqn-mu(j)-02}
\end{eqnarray}
Indeed, for any $j \in \bbS$, there exists some infinite $\bbH_j \subseteq \bbH$ such that
\begin{align}
\lim_{\scriptstyle n\to\infty \atop \scriptstyle n \in \bbH_j}
\pi_n(j)
= 
\limsup_{\scriptstyle n\to\infty \atop \scriptstyle n \in \bbH}
\pi_n(j).
\label{eqn-191014-01}
\end{align}
Note that $\bbH_j \subseteq \bbH$ and $\bbH \subseteq \bbH$ (any set is a subset of itself). Therefore, if (\ref{eqn-mu(j)-02}) holds for all infinite $\bbH' \subseteq \bbH$, then
\begin{align}
\liminf_{\scriptstyle n\to\infty \atop \scriptstyle n \in \bbH_j}
\pi_n(j)
=\liminf_{\scriptstyle n\to\infty \atop \scriptstyle n \in \bbH}
\pi_n(j)
= {\alpha \over \pi(j_0)} \pi(j)
\le \pi(j),\qquad j \in \bbS.
\label{eqn-191014-02}
\end{align}
Combining (\ref{eqn-191014-01}) and (\ref{eqn-191014-02}) yields 
\begin{align*}
\limsup_{\scriptstyle n\to\infty \atop \scriptstyle n \in \bbH}
\pi_n(j)
&= \lim_{\scriptstyle n\to\infty \atop \scriptstyle n \in \bbH_j}
\pi_n(j)
= \liminf_{\scriptstyle n\to\infty \atop \scriptstyle n \in \bbH_j}
\pi_n(j)
\nonumber
\\
&
= \liminf_{\scriptstyle n\to\infty \atop \scriptstyle n \in \bbH}
\pi_n(j)
= {\alpha \over \pi(j_0)} \pi(j)
\le \pi(j),\qquad j \in \bbS,
\end{align*}
which shows that (\ref{eqn-lim-pi_n(j)}) and (\ref{ineqn-alpha}) hold. 

To complete the proof of (\ref{eqn-lim-pi_n(j)}) and (\ref{ineqn-alpha}), 
we prove that (\ref{eqn-mu(j)-02}) holds for all infinite $\bbH' \subseteq \bbH$. Let $\bbH'$ be an arbitrary infinite subset of $\bbH$, and let
\begin{equation}
\mu(j)
=\liminf_{\scriptstyle n\to\infty \atop \scriptstyle n \in \bbH'}\pi_n(j),
\qquad j \in \bbS.
\label{defn-mu(j)} 
\end{equation}
From (\ref{defn-alpha}) and (\ref{defn-mu(j)}), we then have
\begin{align}
\alpha 
&=\lim_{\scriptstyle n\to\infty \atop \scriptstyle n \in \bbH}\pi_n(j_0) 
= \lim_{\scriptstyle n\to\infty \atop \scriptstyle n \in \bbH'}\pi_n(j_0)
= \liminf_{\scriptstyle n\to\infty \atop \scriptstyle n \in \bbH'}\pi_n(j_0)
= \mu(j_0).
\label{eqn-alpha}
\end{align}
Furthermore, using (\ref{defn-mu(j)}), $q(j,j)=\lim_{n\to\infty}q_n(j,j)$ (due to Assumption~\ref{assumpt-appendix}~(ii)), $\vc{\pi}_n\vc{Q}_n=\vc{0}$, and Fatou's lemma, we obtain
\begin{eqnarray*}
-\mu(j) q(j,j)
&=&
-\liminf_{\scriptstyle n\to\infty \atop \scriptstyle n \in \bbH'}
\pi_n(j) q_n(j,j) 
= \liminf_{\scriptstyle n\to\infty \atop \scriptstyle n \in \bbH'}
\sum_{i\in\bbS,\, i\neq j} \pi_n(i) q_n(i,j)
\nonumber
\\
&\ge& 
\sum_{i\in\bbS,\, i\neq j} 
\liminf_{\scriptstyle n\to\infty \atop \scriptstyle n \in \bbH'}
\pi_n(i) q_n(i,j)
\nonumber
\\
&=& 
\sum_{i\in\bbS,\, i\neq j} \mu(i) q(i,j),
\qquad j \in \bbS,
\end{eqnarray*}
which leads to
\[
\sum_{i\in\bbS} \mu(i) q(i,j) \le 0, \qquad j \in \bbS.
\]
It thus follows from Proposition~\ref{prop-subinvariant-measure}~(iii) and (iv) that if $\vc{\mu}:=(\mu(j))_{j\in\bbS} \neq \vc{0}$ then $\vc{\mu}$ is a unique (up to constant multiples) invariant measure of the ergodic $Q$-matrix $\vc{Q}$. Therefore, unifying both cases that $\vc{\mu} \neq \vc{0}$ and $\vc{\mu} = \vc{0}$, we can see that there exists some $c \ge 0$ such that
\begin{eqnarray}
\mu(j)
= c \pi(j) \quad \mbox{for all $j \in \bbS$}.
\label{eqn-mu(j)}
\end{eqnarray}
Substituting (\ref{eqn-mu(j)}) into (\ref{eqn-alpha}) yields
\[
c = {\alpha \over \pi(j_0)}.
\]
Combining this with (\ref{defn-mu(j)}) and (\ref{eqn-mu(j)}) results in
\begin{align}
\liminf_{\scriptstyle n\to\infty \atop \scriptstyle n \in \bbH'}\pi_n(j)
&= \mu(j) 
= c\pi(j)
= {\alpha \over \pi(j_0)}\pi(j), \qquad j \in \bbS.
\label{eqn-191014-03}
\end{align}
In addition, it follows from $\vc{\pi}_n\vc{e}=\vc{\pi}\vc{e}=1$, Fatou's lemma, and (\ref{eqn-191014-03}) that
\begin{eqnarray*}
1 = \liminf_{\scriptstyle n\to\infty \atop \scriptstyle n \in \bbH'}
\sum_{j\in\bbS} \pi_n(j)
&\ge& \sum_{j\in\bbS} 
\liminf_{\scriptstyle n\to\infty \atop \scriptstyle n \in \bbH'}\pi_n(j)
={\alpha \over \pi(j_0)} \sum_{j\in\bbS} \pi(j) 
= {\alpha \over \pi(j_0)}.
\end{eqnarray*}
This result and (\ref{eqn-191014-03}) imply (\ref{eqn-mu(j)-02}). Consequently, (\ref{eqn-lim-pi_n(j)}) and (\ref{ineqn-alpha}) have been proved.

We move on to the proof of the statement (ii). Suppose $\alpha > 0$, and let $\bbA$ be a nonempty and finite $\bbA \subset \bbS$. Since $\bbA$ is finite, it follows from (\ref{eqn-lim-pi_n(j)}) that 
\begin{align*}
\lim_{\scriptstyle n\to\infty \atop \scriptstyle n \in \bbH}
\vc{\pi}_n^{\bbA}
&= {\alpha \over \pi(j_0)} \vc{\pi}^{\bbA},
\\
\lim_{\scriptstyle n\to\infty \atop \scriptstyle n \in \bbH}
\vc{\pi}_n^{\bbA}\vc{e}
&= {\alpha \over \pi(j_0)} \vc{\pi}^{\bbA}\vc{e} > 0,
\end{align*}
where the positivity of the second limit is due to $\alpha > 0$, $\vc{\pi}>\vc{0}$, and $\bbA \neq \varnothing$. The above two equations yield
\begin{align*}
\lim_{\scriptstyle n\to\infty \atop \scriptstyle n \in \bbH}
{\vc{\pi}_n^{\bbA} \over \vc{\pi}_n^{\bbA} \vc{e} }
=
{\vc{\pi}^{\bbA} \over \vc{\pi}^{\bbA} \vc{e} },
\end{align*}
which implies that (\ref{eqn-lim-pi_n^A}) holds. The proof has been completed.
\end{proof}

\medskip

The following theorem presents a necessary and sufficient condition that the sequence $\{\vc{\pi}_n\}$ of approximations converges to the stationary distribution vector $\vc{\pi}$. 
\begin{THM}\label{thm-appen-03}
Under Assumption~\ref{assumpt-appendix}, $\lim_{n\to\infty}\|\vc{\pi}_n - \vc{\pi} \|_1 = 0$ if and only if 
\begin{equation}
\liminf_{n\to\infty} 
\sum_{j\in\bbA} \pi_n(j)
\ge \sum_{j\in\bbA} \pi(j)\quad\mbox{for some finite and nonempty $\bbA \subset \bbS$}.
\label{eqn-liminf-sum-pi_n}
\end{equation}
\end{THM}

\begin{proof}
We note that if $\lim_{n\to\infty}\|\vc{\pi}_n - \vc{\pi} \|_1 = 0$ holds then
\begin{equation}
\liminf_{n\to\infty}\pi_n(j) \ge \pi(j)
\quad \mbox{for all $j \in \bbS$},
\label{eqn-liminf-pi_n(j)}
\end{equation}
which implies (\ref{eqn-liminf-sum-pi_n}). Therefore, to complete the proof, we  prove that ``(\ref{eqn-liminf-sum-pi_n}) $\Rightarrow$ (\ref{eqn-liminf-pi_n(j)})", and then prove that ``(\ref{eqn-liminf-pi_n(j)}) $\Rightarrow \lim_{n\to\infty}\|\vc{\pi}_n - \vc{\pi} \|_1 = 0$".

We prove that ``(\ref{eqn-liminf-sum-pi_n}) $\Rightarrow$ (\ref{eqn-liminf-pi_n(j)})". To this end, suppose that (\ref{eqn-liminf-sum-pi_n}) holds while (\ref{eqn-liminf-pi_n(j)}) does not hold, that is,
\begin{equation}
\alpha:=\liminf_{n\to\infty}\pi_n(j_0) < \pi(j_0)
\quad \mbox{for some $j_0 \in \bbS$}.
\label{ineqn-c<=1}
\end{equation}
Let $\bbH$ be an arbitrary infinite subset $\bbH \subset \bbZ_+$ such that
\begin{align*}
\lim_{\scriptstyle n\to\infty \atop \scriptstyle n \in \bbH}\pi_n(j_0) 
&= \liminf_{n\to\infty}\pi_n(j_0) = \alpha < \pi(j_0). 
\end{align*}
It then follows from Lemma~\ref{lem-appen-02}~(i) that
\begin{equation}
\lim_{\scriptstyle n\to\infty \atop \scriptstyle n \in \bbH}
\pi_n(j) 
= {\alpha \over \pi(j_0)}\pi(j) 
< \pi(j) \quad \mbox{for all $j \in \bbS$}.
\label{eqn-lim-pi_n(j)<pi(j)}
\end{equation}
It also follows from (\ref{eqn-lim-pi_n(j)<pi(j)}) and the finiteness of $\bbA$, and $\bbH \subset \bbZ_+$ that
\begin{eqnarray*}
\sum_{j\in\bbA}\pi(j)
> \lim_{\scriptstyle n\to\infty \atop \scriptstyle n \in \bbH}
\sum_{j\in\bbA} \pi_n(j)
&=& \liminf_{\scriptstyle n\to\infty \atop \scriptstyle n \in \bbH}
\sum_{j\in\bbA} \pi_n(j)
\nonumber
\\
&\ge& \liminf_{\scriptstyle n\to\infty \atop \scriptstyle n \in \bbZ_+}
\sum_{j\in\bbA} \pi_n(j)
= \liminf_{n\to\infty} \sum_{j\in\bbA} \pi_n(j),
\end{eqnarray*}
which contradicts (\ref{eqn-liminf-sum-pi_n}). Therefore, (\ref{eqn-liminf-sum-pi_n}) implies (\ref{eqn-liminf-pi_n(j)}).

We prove that ``(\ref{eqn-liminf-pi_n(j)}) $\Rightarrow \lim_{n\to\infty}\|\vc{\pi}_n - \vc{\pi} \|_1 = 0$". For any fixed $j \in \bbS$, there exist some infinite subset $\bbH_j \subset\bbZ_+$ such that
\begin{align*}
\lim_{\scriptstyle n\to\infty \atop \scriptstyle n \in \bbH_j}
\pi_n(j)
&=
\limsup_{n\to\infty}\pi_n(j).
\end{align*}
Thus, applying Lemma~\ref{lem-appen-02}~(i) to the sequence $\{\vc{\pi}_n;n\in\bbH_j\}$ yields
\begin{align*}
\pi(j)
\ge \lim_{\scriptstyle n\to\infty \atop \scriptstyle n \in \bbH_j}
\pi_n(j)
=
\limsup_{n\to\infty}\pi_n(j) 
\quad \mbox{for any fixed $j \in \bbS$}.
\end{align*}
Combining this and (\ref{eqn-liminf-pi_n(j)}) leads to
\begin{align}
\lim_{n\to\infty}\pi_n(j) = \pi(j)\quad \mbox{for all $j \in \bbS$}.
\label{eqn:22/0102-01}
\end{align}
Note here that $\vc{\pi}_n$ and $\vc{\pi}$ are probability vectors. Therefore, due to the dominated convergence theorem, (\ref{eqn:22/0102-01}) implies that $\lim_{n\to\infty}\|\vc{\pi}_n - \vc{\pi} \|_1 = 0$. It has been proved that ``(\ref{eqn-liminf-pi_n(j)}) $\Rightarrow \lim_{n\to\infty}\|\vc{\pi}_n - \vc{\pi} \|_1 = 0$", which completes the proof.
\end{proof}

\begin{REM}
Lemma~\ref{lem-appen-02}~(i) is the $Q$-matrix-version of \cite[Lemma~2.1]{Wolf80} for stochastic matrices. Furthermore, Lemma~\ref{lem-appen-02}~(ii) and   Theorem~\ref{thm-appen-03} are the $Q$-matrix-versions of Corollary~2.2~(i) and (ii), respectively, in \cite{Wolf80}.
\end{REM}


\section*{Acknowledgments}
The author thanks Dr. Masakiyo Miyazawa for valuable comments on the presentation of the notion of {\it quasi-algorithmic constructions}. The research of the author was supported in part by JSPS KAKENHI Grant Number JP21K11770.

%
%
%
\bibliographystyle{plain} 


\end{document}